\documentclass[11pt, a4paper,twoside]{article}
\usepackage{hyperref}
\usepackage{color}
\usepackage{enumerate}
\usepackage{amsmath}
\usepackage{fix-cm}
\usepackage{pb-diagram}
\usepackage[english]{babel}
\usepackage{amssymb,latexsym}
\usepackage{latexsym}
\def \bui#1#2{\mathrel{\mathop{\kern 0pt#1}\limits^{#2}}}

\let\<\langle
\let\>\rangle
\newcommand{\R}{{\mathbb R}}
\renewcommand{\Re}{\mathrm{Re}}

\newtheorem{example}{Examples}[section]
\newtheorem{thm}{Theorem}[section]

\newtheorem{prop}[thm]{Proposition}

\newtheorem{cor}[thm]{Corollary}
\newtheorem{remark}[thm]{Remark}
\newtheorem{remarks}[thm]{Remarks}
\newtheorem{definition}[thm]{Definition}
\newtheorem{notation}[thm]{Notation}
\newtheorem{exabout:ample}[thm]{Example}

\title{New eigenvalue estimates involving Bessel functions}
\author{Fida El Chami\footnote{Lebanese University, Faculty of Sciences II, Department of Mathematics, P.O. Box 90656 Fanar-Matn, Lebanon,
E-mail: \texttt{fchami@ul.edu.lb, ghabib@ul.edu.lb}},\, Nicolas Ginoux\footnote{Universit\'e de Lorraine, CNRS, IECL, F-57000 Metz, E-mail: \texttt{nicolas.ginoux@univ-lorraine.fr}},\, Georges Habib\footnotemark[1]
}

\begin{document}
\date{}
\maketitle
\begin{abstract}  
\noindent Given a compact Riemannian manifold $(M^n,g)$ with boundary $\partial M$, we give an estimate for the quotient $\displaystyle\frac{\int_{\partial M} fd\mu_g}{\int_M fd\mu_g}$, where $f$ is a smooth positive function defined on $M$ that satisfies some inequality involving the scalar Laplacian. By the mean value lemma established in \cite{S1}, we provide a differential inequality for $f$ which, under some curvature assumptions, can be interpreted in terms of Bessel functions. As an application of our main result, a direct proof is given of the Faber-Krahn inequalities for Dirichlet and Robin Laplacian. Also, a new estimate is established for the eigenvalues of the Dirac operator that involves a positive root of Bessel function besides the scalar curvature. Independently, we extend the Robin Laplacian on functions to differential forms. We prove that this natural extension defines a self-adjoint and elliptic operator whose spectrum is discrete and consists of positive real eigenvalues. In particular, we characterize its first eigenvalue and provide a lower bound of it in terms of Bessel functions.
\end{abstract}

\noindent\begin{small}{\it Mathematics Subject Classification} (2010): 53C27, 53C21, 58J60, 35P15, 34B09, 33C10.
\end{small}

\noindent\begin{small}{\it Keywords}: Bessel functions, Eigenvalues, Dirac operator, Yamabe operator, Robin Laplacian. 
\end{small}
 
\section{Introduction}\label{s:intro}

Let $(M^n,g)$ be an $n$-dimensional compact Riemannian manifold with nonempty smooth boundary $\partial M$. Denote by $\nu$ the inward unit normal along the boundary and by $\Delta f:=-\mathrm{tr}(\nabla^2 f)$ the Laplace operator applied to a smooth function $f$ on $M$. The aim of this paper is twofold. First, we  study the spectrum of some differential operators that arise naturally on manifolds with boundary and are closely related to the Laplacian. Second, we generalize the Robin eigenvalue problem, which consists in solving $\Delta f=\lambda f$ on $M$ with boundary condition $\partial_\nu f=\tau f$ on $\partial M$ for some fixed positive parameter $\tau,$ to differential forms on the manifold. In particular, we aim at establishing sharp lower bounds for the smallest eigenvalue that depend on new invariants. \\

Our first fundamental result deals with the relationship between the integrals $\int_M fd\mu_g$ and $\int_{\partial M}f d\mu_g$, where $d\mu_g$ denotes the Riemannian density on either $M$ or $\partial M$.
Namely, we prove that, as soon as the Ricci curvature of $M$ is nonnegative and the (inward) mean curvature of $\partial M$ is positive, for any positive smooth function $f$ on $M$ satisfying $\Delta f\leq \lambda f$ on $M$ for some sufficiently small $\lambda>0$, the quotient $\displaystyle{\frac{\int_{\partial M}fd\mu_g}{\int_M fd\mu_g}}$ can be bounded below solely in terms of Bessel functions of $\lambda$ and a lower bound for the mean curvature, see Theorem \ref{thm:bessel}.
Note in particular that no boundary condition is required on $f$ here.\\

This central inequality has numerous important consequences.
Namely, applying Theorem \ref{thm:bessel} to a first eigenfunction of the Laplace operator, we can recover in a straightforward way Faber-Krahn inequalities for the Dirichlet (Corollary \ref{cor:Dirichlet}) and Robin (Corollary \ref{daners}) eigenvalue problems assuming only nonnegative Ricci curvature on $M$ and positive mean curvature along $\partial M$.
The former, that is originally due to G.~Faber \cite{Faber1923} and E.~Krahn \cite{Krahn1925}, can be considered as standard, see e.g. \cite[Thm. 2 p.87]{C}.
The latter, which is an immediate consequence of a new lower bound for the first Robin-eigenvalue (Corollary \ref{robinestimate}), was proved by D.~Daners \cite[Thm. 1.1]{D} for domains in $\R^n$ and by A.~Savo in \cite[Thm. 4]{S} for more general manifolds with suitable curvature bounds.\\

Surprisingly enough, Theorem \ref{thm:bessel} can also be applied to eigenvalue problems for differential operators that can be considered as ``far'' from the scalar Laplacian.
For instance, new lower eigenvalue bounds for the Dirac operator are derived in Theorem \ref{t:estimeigenvDirac} and Corollary \ref{cor:chiestimate} under several boundary conditions.
The former is a direct application of Theorem \ref{thm:bessel} choosing $f$ to be the squared length function of some Dirac-eigenspinor; the latter combines a so-called Hijazi-type estimate due to S.~Raulot \cite[Thm. 1]{R1} with a new lower eigenvalue bound for the Yamabe operator we obtain in Theorem \ref{thm:yamabe}.
All lower bounds for the Dirac operator involve the scalar curvature of $M$ as well as the mean curvature of $\partial M$, as can be expected because of the central role of the Schr\"odinger-Lichnerowicz formula relating the squared Dirac operator with the associated rough Laplacian.
Moreover, our lower bounds enhance and rely on former ones due to D.~Chen \cite[Thm. 3.1]{Ch} for the so-called gAPS boundary condition (generalizing \cite[Thm. 4]{HMZ1} dealing with the APS boundary condition), O.~Hijazi, S.~Montiel and A.~Rold\'an \ for the so-called CHI \cite[Thm. 3]{HMR} as well as the MIT bag \cite[Thm. 4]{HMR} boundary conditions, and that of D.~Chen again \cite[Thm. 3.3]{Ch} for the so-called mgAPS boundary condition (generalizing \cite[Thm. 5]{HMR} dealing with the mAPS boundary condition).
It is worth mentioning here that, while our curvature assumptions are stronger than those required by the authors cited above, the estimates we obtain are also stronger since they allow for nontrivial bounds even in case the scalar curvature of $M$ vanishes at one point.\\

Coming back to the Robin eigenvalue problem, we show that it can be generalized in a natural way to differential forms by requiring the boundary conditions
\begin{equation}\label{eq:robinboundary}
\iota^*(\nu\lrcorner d\omega)=\tau \iota^*\omega \quad\text{and}\quad \iota^*(\nu\lrcorner \omega)=0
\end{equation}
for an eigenform $\omega$ of the Laplace operator $\Delta:=d\delta_g+\delta_g d$ on $p$-forms.
Here and in the following, $\iota\colon\partial M\to M$ denotes the inclusion map.
As for the case of functions, both the so-called absolute and Dirichlet boundary conditions can be seen as particular cases of \eqref{eq:robinboundary}, the former setting $\tau:=0$ and the latter letting $\tau\to\infty$.
We first check in Theorem \ref{p:Robinisellipticandselfadjoint} that, assuming $\tau>0$, the Laplacian $\Delta$ is a self-adjoint and elliptic operator with trivial kernel.
Relying on \cite[Ch. 5, Sec. 9]{TaylorPDE1}, we also give in Proposition \ref{claim:varcharfirsteigenvRobin} a variational characterization of its first eigenvalue $\lambda_{1,p}(\tau)$. 
As a first consequence, we prove in Proposition \ref{pro:neurodiri} that the first eigenvalue of the Robin Laplacian always lies between the corresponding absolute and Dirichlet ones.
With the help of the Bochner formula, we deduce from Theorem \ref{thm:bessel} an estimate for the first eigenvalue $\lambda_{1,p}(\tau)$ in terms of Bessel functions, see Theorem \ref{estimatepform}.
As a by-product, we can also derive a lower bound for the gap between Robin eigenvalues on differential forms (Theorem \ref{t:estimgap}) and a Gallot-Meyer-type estimate in case the curvature operator of $M$ is positive (Theorem \ref{t:GallotMeyerRobin}).\\

The article is organized as follows.
After introducing the necessary preliminaries and notations in Section \ref{s:prelim}, we prove the main inequality \eqref{inequalitybessel} and derive its first consequences in Section \ref{s:laplacianonfunctions}.
Sections \ref{s:eigenvestimDirac} and \ref{s:RobinLapldiffforms} are devoted to the application of the main result to the Dirac and the nonscalar Robin eigenvalue problems respectively.
Definitions and some of the basic properties of Bessel functions are recalled in the appendix.

{\bf Acknowledgment:} The second named author thanks the Agence Universitaire de la Francophonie (AUF), the Lebanese University as well as Ines Kath for their support at different stages of the project.
The third named author is supported by a grant from the DAAD within the program ``Research Stays for University Academics and Scientists'' that he would like to thank.
We also thank Sylvie Paycha, Sara Azzali and Alessandro Savo, for fruitful discussions and interesting remarks or questions on the subject.
{\sl Last but not the least}, we are very grateful to Pedro Freitas for informing us about inequality \eqref{eq:ineqFreitastau0tau1}.

\section{Preliminaries}\label{s:prelim}

Let $(M^n,g)$ be a compact Riemannian manifold of dimension $n$ with nonempty smooth boundary $\partial M$. We denote by $\rho\colon M\to [0,\infty[$ the function defined as $\rho(x):=d(x,\partial M).$ Let us first recall basic properties of this function $\rho$ which mainly are contained in \cite{S1} and are also well explained in \cite{GS, RS}. It is not difficult to check, by the triangle inequality, that the function $\rho$ is Lipschitz and its gradient has unit norm a.e. on $M.$ In general, the function $\rho$ is not of class $C^1$ and therefore its Laplacian does not exist as a smooth function. We denote throughout this paper by $\nu$ the inward unit normal vector field to the boundary.
We denote by ${\rm Cut}(\partial M)$ the so-called cut-locus of $\partial M$ in $M$, which is defined as the set of points in $M$ from which more than one minimizing geodesic to the boundary exists.

This set is closed and has measure zero and moreover the function $\rho$ is smooth on its complement (that is usually called the set of regular points), see e.g. \cite[Thm. D.1]{S1}.
To be more precise, the function $\rho$ is smooth on the set $\rho^{-1}[0,{\rm inj}(\partial M)[$ where ${\rm inj}(\partial M)=d(\partial M, {\rm Cut}(\partial M))$ is the injectivity radius. In this case, it is proved in \cite[Subsection 3.2]{S1} that the Laplacian of $\rho$ splits ``in the distributional sense" into a regular part $\Delta_{\rm reg}\rho$ and a positive singular part (with support in the cut-locus) that both can be computed in terms of the local normal coordinates \cite[Eq. 5]{S1}. More explicitly, if we denote by $(r,x)$  the normal coordinates of any regular point ($r$ being the distance of that regular point to $x\in \partial M$), one has $\Delta_{\rm reg}\rho(r,x)= -\frac{1}{\theta}\frac{\partial\theta}{\partial r}(r,x),$ where $\theta$ denotes the density of the pull-back of the volume form (via the local normal exponential map) in normal coordinates.

Given now any smooth function $f$ on $M,$ we define for any $r\geq 0$ the function $F(r):=\displaystyle\int_{\{\rho> r\}} fd\mu_g.$ Clearly, the function $F$ is Lipschitz and is smooth on the interval $[0,{\rm inj}(\partial M)[.$ Moreover, by the co-area formula, its derivative is given by $F'(r)=-\displaystyle\int_{\{\rho=r\}} f d\mu_g$ a.e. on $[0,\infty[$, see \cite[Lemma 2.4]{S1} for a detailed proof.
The {\it mean value lemma}  expresses the second derivative $F''(r)$ in terms of the Laplacian of $f$ through a differential equation that is valid in the sense of distributions. Namely, \cite[Thm. 2.5]{S1} 
\begin{equation}\label{diffeq}
F''(r)=-\left(\int_{\{\rho>r\}}\Delta f d\mu_g\right)+\rho_*(f\Delta \rho)(r),
\end{equation}
where  $\rho_*(f\Delta \rho)$ denotes the push-forward of $f\Delta\rho$ by $\rho,$ that is for any test-function $\psi$ on $[0,\infty[$ we have
\begin{equation} \label{eqpushforward}
(\rho_*(f\Delta \rho),\psi):=\int_0^\infty\psi(r)\left(\int_{\{\rho=r\}}\Delta\rho(r,x) f(x) d\mu_g\right)dr.
\end{equation}

In order to estimate the push-forward in the mean value lemma by some geometric quantities, we require the manifold $M$ to have $(n-1)K$ as a lower bound of the Ricci curvature and $H_0$ as a lower bound of the mean curvature of the boundary. In this case and by the Heintze-Karcher volume inequalities \cite{HK}, one gets in the sense of distributions that 
\begin{equation}\label{eq:laplaciandistance}
\Delta \rho\geq -\frac{\Theta'}{\Theta}\circ\rho,
\end{equation}
where $\Theta$ is the function defined by $\Theta(r)=(s'_K(r)-H_0s_K(r))^{n-1}$ with 
$$s_K(r):=
\left\{
\begin{array}{lll}
\medskip
\frac{1}{\sqrt{K}}\,\sin(r\sqrt{K}) && \text{if}\,\, K>0, \\
\medskip
r && \text{if}\,\, K=0,\\
\medskip
\frac{1}{\sqrt{|K|}}\,\sinh(r\sqrt{|K|}) && \text{if}\,\, K<0.
\end{array}
\right.
$$ 
Therefore if the function $f$ is nonnegative on $M$, then it follows from \eqref{eqpushforward} that (see also \cite[p. 10]{RS})
\begin{equation}\label{eq:equalityball}
\rho_*(f\Delta \rho)\geq \frac{\Theta'}{\Theta} F',
\end{equation}
on the half line. We point out that when $M$ is a geodesic ball in the simply connected manifold $M_K$ of constant curvature $K,$ then equality holds in \eqref{eq:equalityball} (as well in \eqref{eq:laplaciandistance}) for every smooth function $f$. Also, if we let $R:={\rm max}\{d(x,\partial M)|\, x\in M\}$ be the so-called inner radius of $M,$ then the function $\Theta$ is positive on $[0,R[$ and $\Theta(R)=0$ if and only if $M$ is a geodesic ball in $M_K.$ Moreover, denoting by $\bar R$ the first positive zero of the function $r\mapsto s'_K(r)-H_0s_K(r)$, we have $R\leq \bar R$ and equality holds if and only if $M$ is a ball in $M_K$ \cite[Thm. A]{K}.

In this paper, we are interested in studying solutions of the differential equality \eqref{diffeq} in case $f$ is a positive smooth function satisfying $\Delta f\leq \lambda f$ for some $\lambda\geq 0$ (or later a first eigenfunction of the Laplacian). An easy computation using \eqref{diffeq} and \eqref{eq:equalityball} shows that the corresponding differential inequation arising from such a $f$ is (still in the sense of distributions)
\begin{equation}\label{eq:diffine}
F''(r)-\frac{\Theta'}{\Theta} F'(r)+\lambda F(r)\geq 0,
\end{equation}
with $F(0)=\displaystyle\int_M f d\mu_g$ and $F'(0)=-\displaystyle\int_{\partial M} f d\mu_g.$ Keep in mind here that the manifold $M$ is always assumed to have respectively $(n-1)K$ and $H_0$ as lower bounds of the Ricci curvature and mean curvature. It is well-known from the general theory of differential equations (Gr\"onwall Lemma) that the solution $F$ of \eqref{eq:diffine} is always bigger than or equal to that of the corresponding differential equality -- that is, when \eqref{eq:diffine} is an equality -- with the same initial conditions. However, such differential equations cannot be explicitly solved in general as the term in $\Theta$ is hard to control.
A first step in controlling that term was performed by A. Savo and P. Gu\'erini who compute the infimum of $r\mapsto -\frac{\Theta'}{\Theta}(r)$ over all $r$ running in $[0,R[$ (remember that $F'(r)$ is nonpositive). In this case and under some further curvature condition \cite[Eq. 3.1]{GS}, this infimum turns out to be $(n-1)H_0$ and \eqref{eq:diffine} can be reduced to an inequality with constant coefficients whose corresponding differential equality can be explicitly solved. As a consequence, they find a lower bound for the quotient $\displaystyle\frac{\int_{\partial M} fd\mu_g}{\int_M fd\mu_g}$ in terms of $H_0$ and the inner radius $R$ \cite[Thm. 3.1]{GS}. They also deduce well-known sharp estimates for the first eigenvalue of the Dirichlet Laplacian (when $f$ is an eigenfunction) \cite[Cor. 3.2]{GS}, such as McKean \cite{MK} and Li-Yau \cite{LY} inequalities.

A second step was initiated by S.~Raulot and A.~Savo who consider subharmonic functions (i.e. $\lambda=0$). In this particular case, the corresponding differential equation associated to \eqref{eq:diffine} is a linear first order differential equation in $F'$ and the solution can be expressed in terms of $\Theta.$ Namely, they prove that the quotient $\displaystyle\frac{\int_{\partial M} fd\mu_g}{\int_M fd\mu_g}$ is bounded from below by $\frac{1}{\int_0^R\Theta(r)dr}$ \cite[Thm. 10]{RS}. This leads to an estimate for the first eigenvalue of the so-called biharmonic Steklov operator originally introduced in \cite{KS} and \cite{P}.

As we said before, the expression of $\Theta$ which involves sine and hyperbolic sine is difficult to manage, we shall therefore restrict ourselves to the case where $K=0.$ In this case, the term $\frac{\Theta'}{\Theta}$ becomes equal to $\frac{(n-1)H_0}{1-rH_0}.$ Therefore, if we make the change of variable $s=1-rH_0$ and assume moreover that $H_0>0,$ the corresponding differential equality of $\eqref{eq:diffine}$ becomes an equation of Bowman type (see Equation \eqref{diffequbowm} in the appendix) that can be solved in terms of Bessel functions. It turns out that, depending on the dimension of the manifold, we get solutions depending on Bessel functions of first and second kind.

\section{Laplacian on functions}\label{s:laplacianonfunctions}

\noindent In this section, we establish an eigenvalue estimate for the Dirichlet and the Robin Laplacian. As we said before, we express the quotient $\displaystyle\frac{\int_{\partial M} fd\mu_g}{\int_M fd\mu_g}$ in terms of Bessel functions with the help of the mean value lemma. Here $f$ denotes a positive smooth function satisfying some inequality in terms of the Laplacian. Note that no boundary condition is required on $f$ to estimate this quotient. 
In the following, we will denote by $J_\nu$ the Bessel function of order $\nu$ (see the appendix) and by $j_{\nu,k}$ the $k$-th positive zero of $J_\nu$ for $k>0.$

\begin{thm}\label{thm:bessel} Let $(M^n,g)$ be a compact Riemannian manifold with smooth boundary. Assume that the Ricci curvature is nonnegative and the mean curvature is bounded from below by $H_0>0.$ Assume also that there exists a positive smooth function $f$ satisfying $\Delta f\leq \lambda f$ with $\lambda>0$. Then, if  $\frac{\sqrt{\lambda}}{H_0}< j_{\frac{n}{2},1},$  we have that 
\begin{equation}\label{inequalitybessel}
\int_{\partial M} fd\mu_g\geq \sqrt{\lambda} \frac{J_{\frac{n}{2}-1}(\frac{\sqrt{\lambda}}{H_0})}{J_{\frac{n}{2}}(\frac{\sqrt{\lambda}}{H_0})}\int_M f d\mu_g.
\end{equation}
If equality holds in \eqref{inequalitybessel}, then $(M^n,g)$ is isometric to a geodesic ball in $\mathbb{R}^n$ of radius $R=\frac{1}{H_0}$. Conversely, if $M$ is a geodesic ball of radius $\frac{1}{H_0}$ in $\mathbb{R}^n$ and $\Delta f=\lambda f$ for some smooth function $f$ with $ \frac{\sqrt{\lambda}}{H_0}< j_{\frac{n}{2},1},$ then we have equality in \eqref{inequalitybessel}.
\end{thm} 

{\noindent \it Proof.} Since ${\rm Ric}\geq 0,$ we then consider $\Theta(r)=(1-rH_0)^{n-1}.$ Using Inequality \eqref{eq:diffine}, we find that 
\begin{equation}\label{equationeigen}
F''(r)+\frac{(n-1)H_0}{1-rH_0}F'(r)+\lambda F(r)\geq 0.
\end{equation}
Now, we consider the corresponding differential equation $y(r)$ satisfying 
$$y''(r)+\frac{(n-1)H_0}{1-rH_0}y'(r)+\lambda y(r)= 0$$  
with the same initial conditions as $F$. Namely, 
$$F(0)=y(0)=\int_M fd\mu_g\,\, \text{and}\,\, F'(0)=y'(0)=-\int_{\partial M} f d\mu_g.$$ 
By making a change of variable $s=1-rH_0,$ the above differential equation transforms into a Bowman equation, see Equation \eqref{diffequbowm} in the appendix where we let $\gamma:=1, \beta^2:=\frac{\lambda}{H_0^2}$ and $\alpha:=m:=\frac{n}{2}$.
Thus the solution is given by  
$$y(r)=\left\{
\begin{array}{lll}
\medskip
(1-rH_0)^{\frac{n}{2}}\left(AJ_{\frac{n}{2}}\left(\frac{\sqrt{\lambda}}{H_0}(1-rH_0)\right)+BY_{\frac{n}{2}}\left(\frac{\sqrt{\lambda}}{H_0}(1-rH_0)\right)\right)&& \text{if}\,\, n\,\, \text{is even}, \\
(1-rH_0)^{\frac{n}{2}}\left(AJ_{\frac{n}{2}}\left(\frac{\sqrt{\lambda}}{H_0}(1-rH_0)\right)+BJ_{-\frac{n}{2}}\left(\frac{\sqrt{\lambda}}{H_0}(1-rH_0)\right)\right)&& \text{if}\,\, n\,\, \text{is odd}.\\
\end{array}
\right.
$$  
We will first consider the case when $n$ is odd. Taking into account the initial conditions of $y$, the constants $A$ and $B$ must solve the linear system 
$$\left\{\begin{array}{lll}
\medskip
AJ_{\frac{n}{2}}(\frac{\sqrt{\lambda}}{H_0})+BJ_{-\frac{n}{2}}(\frac{\sqrt{\lambda}}{H_0})=\displaystyle\int_M f d\mu_g \\
-A\sqrt{\lambda}J_{\frac{n}{2}-1}(\frac{\sqrt{\lambda}}{H_0})+B\sqrt{\lambda}J_{-\frac{n}{2}+1}(\frac{\sqrt{\lambda}}{H_0})=-\displaystyle\int_{\partial M} f d\mu_g.
\end{array}
\right.$$
For the second equation in the above system, we use the fact that the derivative of the function $r\mapsto (1-rH_0)^{\frac{n}{2}}J_{\frac{n}{2}}(\frac{\sqrt{\lambda}}{H_0}(1-rH_0))$ is equal to $-\sqrt{\lambda}(1-rH_0)^{\frac{n}{2}}J_{\frac{n}{2}-1}(\frac{\sqrt{\lambda}}{H_0}(1-rH_0))$ from the third equation in \eqref{derivativebessel}. Also we use the fourth equation in \eqref{derivativebessel} to compute the derivative of the function $r\mapsto (1-rH_0)^{\frac{n}{2}}J_{-\frac{n}{2}}(\frac{\sqrt{\lambda}}{H_0}(1-rH_0)).$ 
This linear system has clearly a solution $(A,B)$ since from the first equation in \eqref{relationbessel} the determinant of the corresponding matrix is equal to $\frac{2H_0}{\pi}{\rm sin}(\frac{\pi n}{2})\neq 0,$ as $n$ is assumed to be odd. Therefore, we deduce that 
$$A=\frac{\pi}{2H_0 {\rm sin}(\frac{\pi n}{2})}\left(J_{-\frac{n}{2}}\left(\frac{\sqrt{\lambda}}{H_0}\right)\int_{\partial M} fd\mu_g+\sqrt{\lambda}J_{-\frac{n}{2}+1}\left(\frac{\sqrt{\lambda}}{H_0}\right)\int_M fd\mu_g\right),$$
and, 
$$B=\frac{\pi}{2H_0 {\rm sin}(\frac{\pi n}{2})}\left(-J_{\frac{n}{2}}\left(\frac{\sqrt{\lambda}}{H_0}\right)\int_{\partial M} fd\mu_g+\sqrt{\lambda}J_{\frac{n}{2}-1}\left(\frac{\sqrt{\lambda}}{H_0}\right)\int_M fd\mu_g\right).$$ 
Now Inequality \eqref{equationeigen} allows to deduce that, in the sense of distributions, $F(r)\geq y(r).$ Therefore $R\geq R_0$ where $R_0$ is the first positive zero of $y$ (such an $R_0$ exists since $F(R)=0$). As $\Theta$ is positive on $[0,R[,$ two cases may occur for $\Theta(R_0)$, which is nonnegative:\\
 
\noindent $\bullet$ Case where $\Theta(R_0)=0$: In this case, $R_0$ is the unique positive zero of $\Theta$ (recall that $\Theta(r)=(1-rH_0)^{n-1}$). Hence, from \cite[Prop. 14]{RS} (see also \cite[Thm. A]{K}) we deduce that $R_0\geq R$ and thus $R=R_0=\frac{1}{H_0}$. The manifold $M$ is then isometric to the geodesic ball in $\R^n$ of radius $R.$  Now, from the series expansion of the Bessel function in the appendix, we have that the term $(1-rH_0)^{\frac{n}{2}}J_{\frac{n}{2}}(\frac{\sqrt{\lambda}}{H_0}(1-rH_0))$ tends to $0$ and $(1-rH_0)^{\frac{n}{2}}J_{-\frac{n}{2}}(\frac{\sqrt{\lambda}}{H_0}(1-rH_0))$ tends to the constant term of the series $\frac{(\frac{\sqrt{\lambda}}{2H_0})^{-\frac{n}{2}}}{\Gamma(-\frac{n}{2}+1)}\neq 0$ when $r\to \frac{1}{H_0}$. Thus, the fact that $y(R_0)=y(\frac{1}{H_0})=0$ yields $B=0$ and we deduce the equality in \eqref{inequalitybessel}.

\noindent $\bullet$ Case where $\Theta(R_0)>0$: As $y(R_0)=0$, we get that  
$-\frac{A}{B}=\frac{J_{-\frac{n}{2}}\left(\frac{\sqrt{\lambda}}{H_0}(1-R_0H_0)\right)}{J_{\frac{n}{2}}\left(\frac{\sqrt{\lambda}}{H_0}(1-R_0H_0)\right)}.$ We notice here that we can assume that $B\neq 0$ since otherwise we get the equality in \eqref{inequalitybessel}. Also, by assumption the inequalities
\begin{equation}\label{eq:doubleine}
0<\frac{\sqrt{\lambda}}{H_0}\Theta(R_0)^\frac{1}{n-1}=\frac{\sqrt{\lambda}}{H_0}(1-R_0H_0)<\frac{\sqrt{\lambda}}{H_0}< j_{\frac{n}{2},1}
\end{equation}
assure that $J_{\frac{n}{2}}\left(\frac{\sqrt{\lambda}}{H_0}(1-R_0H_0)\right)\neq 0.$ Now, an easy computation that uses Equations \eqref{derivativebessel} and \eqref{relationbessel} in the appendix shows that $(\frac{J_{-\nu}}{J_\nu})'(x)=-2\frac{{\rm sin}\pi \nu}{\pi xJ_\nu^2}$ for all real $\nu.$ In particular, for $\nu=\frac{n}{2},$ the function $x\mapsto \frac{J_{-\frac{n}{2}}}{J_\frac{n}{2}}(x)$ is increasing (resp. decreasing) when $\frac{n-1}{2}$ is odd (resp. even) on $(0,\infty)\setminus\{{\rm zeros\, of\,} J_\frac{n}{2}\}$. But using the expressions of $A$ and $B$, we deduce from \eqref{eq:doubleine} that the following inequality 
 $$\frac{J_{-\frac{n}{2}}\left(\frac{\sqrt{\lambda}}{H_0}\right)\displaystyle\int_{\partial M} fd\mu_g+\sqrt{\lambda}J_{-\frac{n}{2}+1}\left(\frac{\sqrt{\lambda}}{H_0}\right)\displaystyle\int_M fd\mu_g}{J_{\frac{n}{2}}\left(\frac{\sqrt{\lambda}}{H_0}\right)\displaystyle\int_{\partial M} fd\mu_g-\sqrt{\lambda}J_{\frac{n}{2}-1}\left(\frac{\sqrt{\lambda}}{H_0}\right)\int_M fd\mu_g}< \frac{J_{-\frac{n}{2}}\left(\frac{\sqrt{\lambda}}{H_0}\right)}{J_{\frac{n}{2}}\left(\frac{\sqrt{\lambda}}{H_0}\right)},\,\,\,\, ({\rm resp.}\,\, > )$$
holds when $\frac{n-1}{2}$ is odd (resp. even). Taking the common denominator and using again the first equation in \eqref{relationbessel} yields Inequality  \eqref{inequalitybessel}.\\

\noindent The case when $n$ is even is similar to the odd case where $J_{-\frac{n}{2}}$ is replaced by $Y_\frac{n}{2}.$ In this case, the linear system becomes 
$$\left\{\begin{array}{lll}
\medskip
AJ_{\frac{n}{2}}(\frac{\sqrt{\lambda}}{H_0})+BY_{\frac{n}{2}}(\frac{\sqrt{\lambda}}{H_0})=\displaystyle\int_M fd\mu_g \\
-A\sqrt{\lambda}J_{\frac{n}{2}-1}(\frac{\sqrt{\lambda}}{H_0})-B\sqrt{\lambda}Y_{\frac{n}{2}-1}(\frac{\sqrt{\lambda}}{H_0})=-\displaystyle\int_{\partial M} fd\mu_g,
\end{array}
\right.$$
which clearly admits a solution by taking into account the second equation in \eqref{relationbessel}. The constants $A$ and $B$ are then equal to 
$$A=\frac{-H_0}{2\pi}\left(Y_{\frac{n}{2}}\left(\frac{\sqrt{\lambda}}{H_0}\right)\int_{\partial M} fd\mu_g-\sqrt{\lambda}Y_{\frac{n}{2}-1}\left(\frac{\sqrt{\lambda}}{H_0}\right)\int_M fd\mu_g\right),$$
and, 
$$B=\frac{-H_0}{2\pi}\left(-J_{\frac{n}{2}}\left(\frac{\sqrt{\lambda}}{H_0}\right)\int_{\partial M} fd\mu_g+\sqrt{\lambda}J_{\frac{n}{2}-1}\left(\frac{\sqrt{\lambda}}{H_0}\right)\int_M fd\mu_g\right).$$  
Now we proceed as in the odd case, i.e. two cases occur as well. When $\Theta(R_0)=0$, we deduce that $B=0$ as $(1-rH_0)^{\frac{n}{2}}Y_{\frac{n}{2}}(\frac{\sqrt{\lambda}}{H_0}(1-rH_0))$ tends to the constant term of the series $-\frac{1}{\pi}(\frac{n}{2}-1)!(\frac{\sqrt{\lambda}}{2H_0})^{-\frac{n}{2}}\neq 0$. Thus, we get the equality in \eqref{inequalitybessel}. When $\Theta(R_0)>0,$ we have $-\frac{A}{B}=\frac{Y_{\frac{n}{2}}\left(\frac{\sqrt{\lambda}}{H_0}(1-R_0H_0)\right)}{J_{\frac{n}{2}}\left(\frac{\sqrt{\lambda}}{H_0}(1-R_0H_0)\right)}.$ An easy computation that uses Equations \eqref{derivativebessel} and \eqref{relationbessel} in the appendix shows that $(\frac{Y_{\nu}}{J_\nu})'(x)=\frac{2}{\pi xJ_\nu^2}$ which is always positive on $(0,\infty)\setminus\{{\rm zeros\, of\,} J_\nu\}.$ Hence $x\mapsto \frac{Y_{\nu}}{J_\nu}(x)$ is increasing and thus for $\nu=\frac{n}{2}$, we find after using Inequalities \eqref{eq:doubleine} that
$$\frac{Y_{\frac{n}{2}}\left(\frac{\sqrt{\lambda}}{H_0}\right)\displaystyle\int_{\partial M} fd\mu_g-\sqrt{\lambda}Y_{\frac{n}{2}-1}\left(\frac{\sqrt{\lambda}}{H_0}\right)\displaystyle\int_M fd\mu_g}{J_{\frac{n}{2}}\left(\frac{\sqrt{\lambda}}{H_0}\right)\displaystyle\int_{\partial M} fd\mu_g-\sqrt{\lambda}J_{\frac{n}{2}-1}\left(\frac{\sqrt{\lambda}}{H_0}\right)\int_M fd\mu_g}< \frac{Y_{\frac{n}{2}}\left(\frac{\sqrt{\lambda}}{H_0}\right)}{J_{\frac{n}{2}}\left(\frac{\sqrt{\lambda}}{H_0}\right)},$$ 
which leads to the same result as before by the second equation in \eqref{relationbessel}. In the rest of the proof, we discuss the equality case of Inequality \eqref{inequalitybessel}. Assume namely that equality holds in \eqref{inequalitybessel}, then $B=0.$ Next, we prove that $\Theta(R_0)>0$ (recall that $R_0$ is being the first positive root of $y$) cannot occur in this case, so that we are just left with $\Theta(R_0)=0$ which means $M$ is a geodesic ball of radius $\frac{1}{H_0}$ in $\mathbb{R}^n$. Indeed, because $y(R_0)=0,$ we write 
$$0=y(R_0)=A(1-R_0H_0)^{\frac{n}{2}}J_{\frac{n}{2}}(\frac{\sqrt{\lambda}}{H_0}(1-R_0H_0)).$$
Since $J_{\frac{n}{2}}(\frac{\sqrt{\lambda}}{H_0}(1-R_0H_0))\neq 0$ because of \eqref{eq:doubleine}, if we assume by contradiction that $\Theta(R_0)>0$, then the above equality gives that $A=0$. Thus we get $y=0$ which contradicts the fact that $y(0)=\displaystyle\int_M f d\mu_g>0.$ To finish the last part of the equality case, let $M$ be a geodesic ball of radius $\frac{1}{H_0}$ in $\mathbb{R}^n$. As mentioned before, the inequality in \eqref{eq:equalityball}, as well as in \eqref{equationeigen}, are in this case equalities. Therefore, we get that $F(r)=y(r)$ and $R=R_0.$ Because on the ball we have that $\Theta(R)=0,$ we deduce that $B=0.$
\hfill$\square$

\

\begin{remark}\label{remarkroot}\
{\rm
\begin{enumerate} 
\item Recall from \cite{AS,W} that the zeros of $J_{\nu}$ and $J_{\nu+1}$ satisfy $j_{\nu,1}<j_{\nu+1,1}<j_{\nu,2}<j_{\nu+1,2}<\cdots.$
As a consequence, the function $x\mapsto \frac{J_{\frac{n}{2}-1}}{J_{\frac{n}{2}}}(x)$ is positive on the interval $]0, j_{\frac{n}{2}-1,1}[$ and negative on $]j_{\frac{n}{2}-1,1}, j_{\frac{n}{2},1}[$.
In particular, Inequality \eqref{inequalitybessel} does not provide any new information on the interval $]j_{\frac{n}{2}-1,1}, j_{\frac{n}{2},1}[$. 
\item The case where $H_0=0$ was handled in \cite[Thm. 3.1]{GS} but we just add the result for completeness. If $\sqrt{\lambda}R<\frac{\pi}{2},$ then  the corresponding inequality is
$$\int_{\partial M} fd\mu_g\geq \sqrt{\lambda} {\rm cot}(\sqrt{\lambda}R)\int_M fd\mu_g.$$
In this case, one can deduce estimates for the first eigenvalues of the Dirichlet and Robin Laplacian, see \cite[Cor. 3.2]{GS}, \cite[Cor. 3]{S}, \cite{LY}, \cite{MK}. 
\end{enumerate} }
\end{remark}  

An immediate consequence of Theorem \ref{thm:bessel} is when the function $f$ is subharmonic (i.e. $\Delta f\leq 0$). Namely, we have 
\begin{cor}\label{c:ineqHeintzeKarcherRos} Let $(M^n,g)$ be a compact Riemannian manifold with smooth boundary. Assume that the Ricci curvature of $(M^n,g)$ is nonnegative and the mean curvature of $\partial M$ is bounded from below by $H_0>0.$ Let $f$ be any positive and subharmonic function. Then 
\begin{equation}\label{raulotsavo}
\frac{\int_{\partial M} f d\mu_g}{\int_{M} f d\mu_g}\geq nH_0. 
\end{equation}
Equality holds if and only if $M$ is isometric to the geodesic ball $B_{H_0}$ of radius $\frac{1}{H_0}$.
In particular, for $f=1$, one has
\begin{equation}\label{ros}
\frac{{\rm Vol}(\partial M)}{{\rm Vol}(M)}\geq nH_0, 
\end{equation}
where equality holds if and only if $M$ is isometric to $B_{H_0}$.
\end{cor}
{\noindent \it Proof.} By applying Theorem \ref{thm:bessel} to the function $f$, we deduce from Inequality \eqref{inequalitybessel} that for any $\lambda>0$ with $\frac{\sqrt{\lambda}}{H_0}<j_{\frac{n}{2},1}$ we have 
$$\frac{\int_{\partial M} f d\mu_g}{\int_{M} f d\mu_g}\geq \sqrt{\lambda} \frac{J_{\frac{n}{2}-1}(\frac{\sqrt{\lambda}}{H_0})}{J_{\frac{n}{2}}(\frac{\sqrt{\lambda}}{H_0})}=H_0x\frac{J_{\frac{n}{2}-1}(x)}{J_{\frac{n}{2}}(x)},$$
with $x:= \frac{\sqrt{\lambda}}{H_0}.$ Taking the limit as $x\to 0$ and using the fact that for all $\nu\geq 0$ we have $\frac{J_\nu(x)}{J_{\nu+1}(x)}\approx \frac{2(\nu+1)}{x}$ for small $x$ \cite[p. 192]{L}, leads to the result. The equality case follows also directly from Theorem \ref{thm:bessel}.
\hfill$\square$

\noindent We point out that Inequality \eqref{raulotsavo} is weaker than Raulot-Savo's estimate \cite[Thm. 10]{RS} which states that 
$$\frac{\int_{\partial M} f d\mu_g}{\int_{M} f d\mu_g}\geq \frac{1}{\int_0^R\Theta(r)dr}=\frac{nH_0}{1-(1-RH_0)^n}.$$
This is due to the different solution of the differential equation in \cite{RS} which does not involve the Bessel functions. Also, we notice that Inequality \eqref{ros} is also weaker than the estimate in \cite{HK}, \cite{Ros}, known as Heintze-Karcher-Ros,  which is 
$$\int_{\partial M}\frac{1}{H}  d\mu_g\geq n {\rm Vol}(M),$$
where $H$ is the mean curvature.

\noindent Recall now that the quotient of two consecutive Bessel functions is given by the series \cite[p. 498]{W} 
$$\frac{J_{\nu+1}(x)}{J_\nu(x)}=\sum_{k\geq 1} \frac{2x}{j^2_{\nu,k}-x^2},$$  
for any $\nu>-1.$ Hence, for $\nu=\frac{n}{2}-1$, the function $x\mapsto \frac{J_{\frac{n}{2}}(x)}{J_{\frac{n}{2}-1}(x)}$ increases on $\R\setminus\{{\rm zeros\, of\,} J_{\frac{n}{2}-1}\}$. It is also positive on $]0,j_{\frac{n}{2}-1,1}[$ and negative on $]j_{\frac{n}{2}-1,1},j_{\frac{n}{2},1}[$.
Using this, we obtain the following estimate for the first eigenvalue of the Dirichlet Laplacian, known as Faber-Krahn inequality \cite{C}:

\begin{cor} \label{cor:Dirichlet} Let $(M^n,g)$ be a compact Riemannian manifold with smooth boundary. Assume that the Ricci curvature of $(M^n,g)$ is nonnegative and the mean curvature of $\partial M$ is bounded from below by $H_0>0.$
Let $B_{H_0}$ be the geodesic ball of radius $\frac{1}{H_0}$ in the Euclidean space. Then the first eigenvalue of the Dirichlet Laplacian $\lambda^D_1$ satisfies 
$\lambda^D_1\geq \lambda_1^D(B_{H_0})=H_0^2 j^2_{\frac{n}{2}-1,1}.$ Equality is attained if and only if $M$ is isometric to the ball $B_{H_0}.$
\end{cor}  

{\noindent \it Proof.} Let $f$ be a positive eigenfunction of Dirichlet Laplacian associated to the first eigenvalue $\lambda^D_1$. If $\sqrt{\lambda^D_1}<H_0j_{\frac{n}{2}-1,1}<H_0j_{\frac{n}{2},1},$ then we get from Inequality \eqref{inequalitybessel} that 
\begin{equation*}
0=\int_{\partial M} f d\mu_g\geq \sqrt{\lambda^D_1} \frac{J_{\frac{n}{2}-1}\left(\frac{\sqrt{\lambda^D_1}}{H_0}\right)}{J_{\frac{n}{2}}\left(\frac{\sqrt{\lambda^D_1}}{H_0}\right)}\int_M f d\mu_g>0.
\end{equation*}
This leads to a contradiction. If now equality is realized, then we still have $\frac{\sqrt{\lambda^D_1}}{H_0}< j_{\frac{n}{2},1}$ and thus the inequality
$$0=\int_{\partial M} f d\mu_g\geq \sqrt{\lambda^D_1} \frac{J_{\frac{n}{2}-1}\left(\frac{\sqrt{\lambda^D_1}}{H_0}\right)}{J_{\frac{n}{2}}\left(\frac{\sqrt{\lambda^D_1}}{H_0}\right)}\int_M f d\mu_g=0,$$
becomes an equality. Therefore, we deduce the result from the characterization of the equality case in Theorem \ref{thm:bessel}. This ends the proof. 
\hfill$\square$

In the following, we are interested in estimating the eigenvalues of the Robin Laplacian. Recall that this boundary problem is defined as follows: Fix a positive parameter $\tau$ and consider the boundary value problem 
\begin{equation}\label{eq:robinproblem}
\left\{
\begin{array}{lll}
\medskip
\Delta f=\lambda f\, && \text{on }\,\, M, \\
\medskip
\frac{\partial f}{\partial \nu}=\tau f && \text{on}\,\, \partial M,\\
\end{array}
\right.
\end{equation}  
where $\nu$ is the inward normal vector field to the boundary. It is well-known that the eigenvalues of the Robin Laplacian form an increasing sequence $0<\lambda_1(\tau,M)<\lambda_2(\tau,M)\leq\cdots$ (counted with multiplicities) and depend continuously on $\tau.$ 
When $\tau$ tends to zero, the Robin Laplacian reduces to the Neumann Laplacian while it is the Dirichlet Laplacian when $\tau\to \infty$. Using Theorem \ref{thm:bessel}, we will establish an estimate for the first eigenvalue of the Robin Laplacian in terms of the zeros of Bessel functions. First, observe that, for any eigenfunction $f$ associated to an eigenvalue $\lambda$ of the problem \eqref{eq:robinproblem}, we have
\begin{equation}\label{eq:robin}
\lambda\displaystyle\int_M f d\mu_g=\int_M \Delta f d\mu_g=\int_{\partial M}\frac{\partial f}{\partial \nu}d\mu_g=\tau\int_{\partial M} fd\mu_g.
\end{equation}

\noindent Therefore, the quotient $\displaystyle\frac{\int_{\partial M} fd\mu_g}{\int_M fd\mu_g}$ is just equal to $\frac{\lambda}{\tau}$ whenever $\int_M fd\mu_g>0$.
Taking this fact into account,  we get the following  

\begin{cor} \label{cor:estimate} Let $(M^n,g)$ be a compact Riemannian manifold with smooth boundary. Assume that the Ricci curvature of $(M^n,g)$ is nonnegative and the mean curvature of $\partial M$ is bounded from below by $H_0>0.$ If $\frac{\sqrt{\lambda_1(\tau,M)}}{H_0}< j_{\frac{n}{2},1},$ then  
$$\sqrt{\lambda_1(\tau,M)}\geq \tau \frac{J_{\frac{n}{2}-1}\left(\frac{\sqrt{\lambda_1(\tau,M)}}{H_0}\right)}{J_{\frac{n}{2}}\left(\frac{\sqrt{\lambda_1(\tau,M)}}{H_0}\right)}.$$ 
If  $\frac{\sqrt{\lambda_1(\tau,M)}}{H_0}< j_{\frac{n}{2}-1,1},$ equality is realized if and only if $M$ is isometric to the ball $B_{H_0}.$  
\end{cor} 

\begin{remark}\label{rem:eigenvalueball}
{\rm From the characterization of the equality case of Inequality \eqref{inequalitybessel}, for which $B=0$ as was shown in the proof of Theorem \ref{thm:bessel}, and in view of \eqref{eq:robin}, we deduce that on a geodesic ball $B_{H_0}$ in $\mathbb{R}^n$ any eigenvalue $\lambda(\tau,B_{H_0})$ of the Robin Laplacian associated with an eigenfunction $f$ satisfies the equality 
$$\left[J_{\frac{n}{2}}\left(\frac{\sqrt{\lambda(\tau,B_{H_0})}}{H_0}\right) \sqrt{\lambda(\tau,B_{H_0})}-\tau J_{\frac{n}{2}-1}\left(\frac{\sqrt{\lambda(\tau,B_{H_0})}}{H_0}\right)\right]\int_{B_{H_0}} f d\mu_g=0.$$ 
Hence, for the first positive eigenvalue $\lambda_1(\tau,B_{H_0})$ (in this case $f$ is positive), the term $\frac{\sqrt{\lambda_1(\tau,B_{H_0})}}{H_0}$ is a root of the function $x\mapsto \frac{J_{\frac{n}{2}}(x)}{J_{\frac{n}{2}-1}(x)}-\frac{\tau}{H_0 x}$ which is defined on $\R\setminus\{{\rm zeros\, of\,} J_{\frac{n}{2}-1}\}$ and increases from $-\infty$ to $\infty$ on $]0,j_{\frac{n}{2}-1,1}[$. It is indeed the first positive zero on that interval, see \cite[Rem. 2.9]{AFK1}, \cite[p. 4]{AFK2} for more details. Therefore, we deduce that $\frac{\sqrt{\lambda_1(\tau,B_{H_0})}}{H_0}$ lies in the interval $]0,j_{\frac{n}{2}-1,1}[.$ 
}
\end{remark}

\noindent Using the above corollary, we have the following estimate

\begin{cor} \label{robinestimate} Let $(M^n,g)$ be a compact Riemannian manifold with smooth boundary. Assume that the Ricci curvature of $(M^n,g)$ is nonnegative and the mean curvature of $\partial M$ is bounded from below by $H_0>0.$ Fix any positive number $\tau_0<j_{\frac{n}{2}-1,1}$ and set $\alpha=\sum_{k\geq 1}\frac{2\tau_0^2}{j_{\frac{n}{2}-1,k}^2-\tau_0^2}.$ If $\tau\geq \alpha H_0,$ then 
$$\lambda_1(\tau,M)\geq H_0^2\tau_0^2.$$  
Equality case is realized if and only if $M$ is isometric to $B_{H_0}$ and $\tau=\alpha H_0.$ 
\end{cor}

{\noindent \it Proof.} Assume that $\sqrt{\lambda_1(\tau,M)}< H_0\tau_0<H_0j_{\frac{n}{2}-1,1},$ then by Corollary \ref{cor:estimate} we get that
$$\sqrt{\lambda_1(\tau,M)}\geq \tau \frac{J_{\frac{n}{2}-1}\left(\frac{\sqrt{\lambda_1(\tau,M)}}{H_0}\right)}{J_{\frac{n}{2}}\left(\frac{\sqrt{\lambda_1(\tau,M)}}{H_0}\right)}> \tau \frac{\sqrt{\lambda_1(\tau,M)}}{\alpha H_0},$$ 
with $\displaystyle{\alpha=\tau_0\frac{J_{\frac{n}{2}}(\tau_0)}{J_{\frac{n}{2}-1}(\tau_0)}=\sum_{k\geq 1}\frac{2\tau_0^2}{j_{\frac{n}{2}-1,k}^2-\tau_0^2}}$.
The last inequality comes from the fact that the function $x\mapsto x\frac{J_{\nu+1}}{J_\nu}(x)=\sum_{k\geq 1} \frac{2x^2}{j^2_{\nu,k}-x^2}$ is increasing on $]0,j_{\frac{n}{2}-1,1}[$ and thus for $\nu=\frac{n}{2}-1$ and $x\in\,]0,\tau_0[$, we have $\frac{J_{\frac{n}{2}}}{J_{\frac{n}{2}-1}}(x)<\frac{\alpha}{x}$. This leads to a contradiction. Assume now that the equality case is attained. Then we still have that $\frac{\sqrt{\lambda_1(\tau,M)}}{H_0}<j_{\frac{n}{2}-1,1}$ and the inequality in Corollary \ref{cor:estimate}
$$H_0\tau_0=\sqrt{\lambda_1(\tau,M)}\geq \tau \frac{J_{\frac{n}{2}-1}\left(\frac{\sqrt{\lambda_1(\tau,M)}}{H_0}\right)}{J_{\frac{n}{2}}\left(\frac{\sqrt{\lambda_1(\tau,M)}}{H_0}\right)}=\tau \frac{J_{\frac{n}{2}-1}(\tau_0)}{J_{\frac{n}{2}}(\tau_0)}=\tau\frac{\tau_0}{\alpha}\geq H_0\tau_0,$$
becomes an equality. Therefore, we deduce that $M$ is isometric to $B_{H_0}$ and $\tau=\alpha H_0.$ Conversely, on the geodesic ball $B_{H_0}$ we have equality in the estimate of Corollary \ref{cor:estimate}. Hence we write, for $\tau=\alpha H_0,$ 
$$
\sqrt{\lambda_1(\tau,B_{H_0})}=\tau \frac{J_{\frac{n}{2}-1}\left(\frac{\sqrt{\lambda_1(\tau,B_{H_0})}}{H_0}\right)}{J_{\frac{n}{2}}\left(\frac{\sqrt{\lambda_1(\tau,B_{H_0})}}{H_0}\right)}\leq \alpha H_0\frac{J_{\frac{n}{2}-1}(\tau_0)}{J_{\frac{n}{2}}(\tau_0)}=\tau_0H_0.
$$
Here, we use the fact that the function $x\mapsto \frac{J_{\nu}}{J_{\nu+1}}(x)$ is decreasing. Hence, we get the other side of the estimate and thus the equality is attained.
\hfill$\square$  

\noindent Using the previous result, we will compare the first eigenvalue of the Robin Laplacian on $M$ to the one on the ball $B_{H_0}$. This is known as the Faber-Krahn inequality proved by Daners \cite{D} for Euclidean domains. We have

\begin{cor} \label{daners} Let $(M^n,g)$ be a compact Riemannian manifold with smooth boundary. Assume that the Ricci curvature of $(M^n,g)$ is nonnegative and the mean curvature of $\partial M$ is bounded from below by $H_0>0.$ Let $B_{H_0}$ be the geodesic ball of mean curvature $H_0.$ Then  $$\lambda_1(\tau,M)\geq \lambda_1(\tau,B_{H_0}).$$ 
Equality is realized if and only if $M$ is isometric to the geodesic ball $B_{H_0}.$
\end{cor} 
{\noindent \it Proof.} In view of Remark \ref{rem:eigenvalueball}, we have that $\frac{\sqrt{\lambda_1(\tau,B_{H_0})}}{H_0}<j_{\frac{n}{2}-1,1}$. Therefore, we set $\tau_0:=\frac{\sqrt{\lambda_1(\tau,B_{H_0})}}{H_0}$ in Corollary \ref{robinestimate}. In this case, we get that
$$\alpha=\tau_0\frac{J_{\frac{n}{2}}(\tau_0)}{J_{\frac{n}{2}-1}(\tau_0)}=\frac{\sqrt{\lambda_1(\tau,B_{H_0})}}{H_0}\frac{\tau}{\sqrt{\lambda_1(\tau,B_{H_0})}}=\frac{\tau}{H_0}.$$
Hence Corollary \ref{robinestimate} finishes the proof of the result.
\hfill$\square$

\begin{remark} \label{remark:lowerbound} \
{\rm
According to \cite{QW}, the best possible lower bound for $j_{\nu,1}$ (for $\nu>0$) is the positive number $\tau_0=\nu-\frac{a_1}{2^\frac{1}{3}}\nu^\frac{1}{3}$ where $a_1\simeq -2.3381$ is the first negative zero of the Airy function. Therefore for $\nu=\frac{n}{2}-1$ with $n\geq 3$, one can easily check that $\tau_0> \frac{n}{2}$. Thus if we choose $\tau\geq \alpha H_0,$ one gets 
$$\lambda_1(\tau,M)\geq \tau_0^2 H_0^2> \frac{n^2}{4}H_0^2> nH_0\tau-\tau^2> (n-1)H_0\tau-\tau^2.$$ 
The last lower bound has been obtained in \cite{BMMP} under the further assumption that $II+\tau>0,$ where $II$ denotes the second fundamental form of the boundary. }
\end{remark}

\section{Eigenvalue estimates for the Dirac operator}\label{s:eigenvestimDirac}

\noindent In this section we give, under curvature assumptions, new estimates for the first eigenvalue of the Dirac operator defined on compact manifolds with boundary. These estimates are expressed in terms of zeros of Bessel functions and a lower bound of the scalar curvature. They improve Friedrich-type estimates originally established on closed manifolds \cite{F} and generalized later on manifolds with boundary, see e.g. \cite{HMR} or \cite[Ch. 4]{GinsurveyDirac} for references.

\noindent We assume here the smooth compact Riemannian manifold $(M^n,g)$ to be spin with fixed spin structure.
For more details on spin manifolds, we refer to e.g. \cite{BFGK91}, \cite{Friedrichbook2000}, \cite[Ch. 1\&2]{LawsonMichelsohn89}, \cite[Ch. 1]{BHMM}, \cite[Ch. 1]{GinsurveyDirac}.
Under that assumption, there exists a smooth Hermitian vector bundle $\Sigma M\longrightarrow M$ called the spinor bundle of $M$ on which $TM$ acts by Clifford multiplication.
We denote by $X\otimes\psi\mapsto X\cdot \psi$ that Clifford multiplication.
There exists on $\Sigma M$ a metric connection that preserves the Clifford multiplication and that we denote by $\nabla^{\Sigma M}$.
The Clifford trace of that connection is the first-order differential operator called the Dirac operator.
Formally, for any section $\psi$ of $\Sigma M$, we have $D\psi=\sum_{j=1}^ne_j\cdot\nabla_{e_j}^{\Sigma M}\psi$, where $(e_j)_{1\leq j\leq n}$ is an arbitrary local g-o.n.b. of $TM$.
 Recall also that a spin structure on $M$ induces a spin structure on $\partial M$ via the inner unit normal vector field $\nu$ along $\partial M$.
This provides a unitary isomorphism
\[\Sigma M_{|_{\partial M}}\longrightarrow\left\{\begin{array}{ll}\Sigma \partial M&\textrm{ if $n$ is odd}\\\Sigma \partial M\oplus\Sigma\partial M&\textrm{ if $n$ is even}\end{array}\right.\]
for which a Gau\ss{}-type-formula relates the compatible connections on $\Sigma M$ and $\Sigma\partial M$.
In particular, the following formula holds along $\partial M$ for any $\psi\in\Gamma(\Sigma M)$, see e.g. \cite[Eq. (1.22)]{GinsurveyDirac}:
\begin{equation}\label{eq:DiracMdM}
D\psi=\nu\cdot\nabla_{\nu}^{\Sigma M}\psi+\nu\cdot\left(D^{\partial M}\psi-\frac{(n-1)H}{2}\psi\right),
\end{equation}
where $H:=\frac{1}{n-1}\mathrm{tr}(II)$ is the mean curvature of $\partial M$ in $M$ and $D^{\partial M}$ is either the Dirac operator on $\partial M$ (if $n$ is odd) or its symmetrization (if $n$ is even), we refer to e.g. \cite[Sec. 1.4]{GinsurveyDirac} for details.
Here $II:=-\nabla^M\nu$ denotes the Weingarten map of the boundary.

\noindent The Dirac operator is known to admit the following four elliptic boundary conditions: CHI, MIT bag, gAPS and mgAPS, see e.g. \cite{HMR} or \cite[Sec. 1.5]{GinsurveyDirac} for a survey. Recall first that the corresponding boundary value problem consists in solving $D\psi=\lambda\psi$ on $M$ where $\psi$ lies in the kernel of the boundary operator $B$ corresponding to one of the boundary condition listed above.
It is proved in \cite[Prop. 1]{HMR} that, under these boundary conditions, the spectrum of the Dirac operator consists of a discrete unbounded sequence of eigenvalues with finite dimensional eigenspaces. For the CHI, gAPS and mgAPS boundary conditions, the spectrum is real, however for the MIT bag condition it is contained in the upper half of the complex line $\mathbb{C}.$

\noindent Let us now recall briefly these boundary conditions. The CHI boundary condition is associated to the so-called {\it chirality operator}, defined by the endomorphism $B_{\rm CHI}:=\frac{1}{2}({\rm Id}-\nu\cdot\mathcal{G}),$ where $\nu$ is the unit normal vector field to $\partial M$ and $\mathcal{G}$ is the restriction on $\partial M$ of the endomorphism $\mathcal{G}:\Sigma M\to \Sigma M$ which is involutive, unitary, parallel and anticommuting with Clifford multiplication on $M$ (it corresponds to the complex volume form for $n$ even). The MIT bag condition is defined by the operator $B_{\rm MIT}:=\frac{1}{2}({\rm Id}-i\nu\cdot).$ For the gAPS, known as {\it generalized Atiyah-Patodi Singer}, the boundary operator $B_{\rm gAPS}$ is defined as the $L^2$-orthogonal projection onto the subspace generated by the eigenvectors of the Dirac operator on $\partial M$ (if $n$ is even or its symmetrization if $n$ is odd) corresponding with eigenvalues not smaller than some number $\beta\leq 0.$ Finally, the boundary operator $B_{\rm mgAPS}$ for the condition, known as {\it modified generalized Atiyah-Patodi-Singer}, is defined as $B_{\rm mgAPS}:=B_{\rm gAPS}({\rm Id}+\nu\cdot).$

\noindent  In \cite{HMR}, the authors provide a Friedrich-type lower bound involving scalar curvature \cite{F} for the first eigenvalue of the Dirac operator and for each of the above boundary conditions (see also \cite{Ch}, \cite{HMZ1}). They also discuss the equality case of those estimates which turns out not to be always achieved depending on the imposed boundary condition; moreover, for the cases where the equality is realized, the boundary has to be minimal. We notice here that the positivity of the scalar curvature as well as the nonnegativity of the mean curvature of the boundary are required in this context in order for the lower bound to be positive. In the following, we will give a new estimate for the eigenvalues of the Dirac operator under the boundary conditions mentioned above (except the MIT bag) by using the result in Theorem \ref{thm:bessel}.
The new fact in our estimate is that the lower bound not only depends on the minimum of the scalar curvature (as for Friedrich's lower bound) but also on a positive root of some function involving Bessel functions.
In particular, this estimate still gives us information on the spectrum when the scalar curvature of the manifold vanishes at one point. 

\begin{thm}\label{t:estimeigenvDirac}
Let $(M^n,g)$ be any smooth compact Riemannian spin manifold with nonempty boundary $\partial M$.
Assume $\mathrm{Ric}\geq0$ on $M$ and $H\geq H_0$ on $\partial M$ for some positive constant $H_0$.
Let $\lambda$ be any eigenvalue of the Dirac operator of $M$ endowed with any of the CHI, gAPS or mgAPS boundary conditions.
Let $\tau_0$ be the only zero of $x\mapsto x\frac{J_{\frac{n}{2}(x)}}{J_{\frac{n}{2}-1}(x)}-(n-1)$ on $]0,j_{\frac{n}{2}-1,1}[$.
Then 
\begin{equation}\label{eq:estimeigenvDirac}
\lambda^2>\frac{n}{4(n-1)}\min_M(S)+\frac{n H_0^2\tau_0^2}{2(n-1)},
\end{equation}
where $S$ is the scalar curvature of $(M^n,g)$.
\end{thm}

\noindent {\it Proof.} We consider the nonnegative function $f:=\frac{1}{2}|\psi|^2$ on $M$, where $\psi$ is any spinor field on $\Sigma M.$ It is elementary to show that, with the help of the Schr\"odinger-Lichnerowicz formula $D^2=(\nabla^{\Sigma M})^*\nabla^{\Sigma M}+\frac{S}{4}\mathrm{Id}$, the following identity holds for $\psi$:
\begin{equation}\label{Diracschro}
\Delta f=-\frac{S}{2}f+\Re(\langle D^2\psi,\psi\rangle)-\frac{1}{n}|D\psi|^2-|P\psi|^2,
\end{equation}
where $P\colon\Gamma(\Sigma M)\longrightarrow\Gamma(T^*M\otimes\Sigma M)$ is the so-called Penrose operator, defined by $P\psi=\nabla_X^{\Sigma M}\psi+\frac{1}{n}X\cdot D\psi$ for any vector field $X\in TM.$
Taking $\psi$ to be a nonzero Dirac-eigenspinor associated to the eigenvalue $\lambda$ (remember that $\lambda$ is real under the imposed boundary conditions), we obtain 
\begin{eqnarray}\label{laplacianspinor}
\Delta f&=&-\frac{S}{2}f+\lambda^2|\psi|^2-\frac{\lambda^2}{n}|\psi|^2-|P\psi|^2\nonumber\\
&=&-\frac{S}{2}f+2\lambda^2f-\frac{2\lambda^2}{n}f-|P\psi|^2\nonumber\\
&=&\frac{2(n-1)}{n}\left(\lambda^2-\frac{n}{4(n-1)}S\right)f-|P\psi|^2\nonumber\\
&\leq&\frac{2(n-1)}{n}\left(\lambda^2-\frac{n}{4(n-1)}\min_M(S)\right)f,
\end{eqnarray}
that is, $\Delta f\leq\mu f$ where $\mu:=\frac{2(n-1)}{n}\left(\lambda^2-\frac{n}{4(n-1)}\min_M(S)\right)$.
Notice that, by all eigenvalue estimates proved in \cite{HMR} for the boundary conditions assumed in our theorem, we have $\mu>0$ (remember that $H\geq H_0>0$). On the other hand, using Gau\ss{} formula \eqref{eq:DiracMdM}, we compute
\begin{eqnarray*}
\int_M\Delta f d\mu_g&=&\int_{\partial M}\partial_\nu fd\mu_g\\
&=&\int_{\partial M}\Re(\langle\nabla_\nu^{\Sigma M}\psi,\psi\rangle)d\mu_g\\
&\bui{=}{\rm\eqref{eq:DiracMdM}}&\int_{\partial M}\Re(\langle-\nu\cdot D\psi-D^{\partial M}\psi+\frac{(n-1)H}{2}\psi,\psi\rangle)d\mu_g\\
&=&-\lambda\int_{\partial M}\underbrace{\Re(\langle\nu\cdot\psi,\psi\rangle)}_{0}d\mu_g-\int_{\partial M}\langle D^{\partial M}\psi,\psi\rangle d\mu_g+\frac{n-1}{2}\int_{\partial M}H|\psi|^2d\mu_g\\
&=&-\int_{\partial M}\langle D^{\partial M}\psi,\psi\rangle d\mu_g+(n-1)\int_{\partial M}Hfd\mu_g.
\end{eqnarray*}
Now since by assumption $H\geq H_0$, $f\geq0$ along $\partial M$ and $\displaystyle\int_{\partial M}\langle D^{\partial M}\psi,\psi\rangle d\mu_g\leq0$ for any of the boundary conditions under consideration (see e.g. \cite[Ch. 4]{GinsurveyDirac}), we obtain
\begin{equation}\label{eq:ineqDeltafDirac}
\int_M\Delta f d\mu_g\geq(n-1)H_0\int_{\partial M}f d\mu_g.
\end{equation}
Notice here that no condition on the Ricci curvature is required to get Inequalities \eqref {laplacianspinor} and \eqref{eq:ineqDeltafDirac}. By contradiction let us now assume that $\frac{\sqrt{\mu}}{H_0}<\tau_0$, where $\tau_0$ is the only zero of $x\mapsto x\frac{J_{\frac{n}{2}(x)}}{J_{\frac{n}{2}-1}(x)}-(n-1)$ on $]0,j_{\frac{n}{2}-1,1}[.$ 
Since in particular $\frac{\sqrt{\mu}}{H_0}<j_{\frac{n}{2}-1,1}$ and by assumption $\mathrm{Ric}\geq0$ on $M$ and $H\geq H_0>0$ on $\partial M$, Theorem \ref{thm:bessel} can be applied to $f$ and yields 
\[\int_{\partial M}fd\mu_g\geq\sqrt{\mu}\frac{J_{\frac{n}{2}-1}(\frac{\sqrt{\mu}}{H_0})}{J_{\frac{n}{2}}(\frac{\sqrt{\mu}}{H_0})}\int_M fd\mu_g.\]
Note that in particular the function $f$ cannot vanish identically on the boundary.
But \eqref{eq:ineqDeltafDirac} together with $\Delta f\leq \mu f$ implies $\mu\displaystyle\int_M fd\mu_g\geq(n-1)H_0\int_{\partial M} f d\mu_g$, so that
\[\int_{\partial M}fd\mu_g\geq\sqrt{\mu}\frac{J_{\frac{n}{2}-1}(\frac{\sqrt{\mu}}{H_0})}{J_{\frac{n}{2}}(\frac{\sqrt{\mu}}{H_0})}\cdot\frac{(n-1)H_0}{\mu}\int_{\partial M}fd\mu_g,\]
from which $\frac{1}{n-1}\geq\frac{J_{\frac{n}{2}-1}(\frac{\sqrt{\mu}}{H_0})}{\frac{\sqrt{\mu}}{H_0}J_{\frac{n}{2}}(\frac{\sqrt{\mu}}{H_0})}$ follows.
By assumption on $\mu$ and since $x\mapsto \frac{J_{\frac{n}{2}-1}(x)}{xJ_{\frac{n}{2}}(x)}$ is decreasing on $]0,j_{\frac{n}{2}-1,1}[$, we deduce that $\frac{1}{n-1}>\frac{J_{\frac{n}{2}-1}(\tau_0)}{\tau_0 J_{\frac{n}{2}}(\tau_0)}=\frac{1}{n-1}$, which is a contradiction.
Therefore $\sqrt{\mu}\geq H_0\tau_0$, which concludes the proof of the inequality \eqref{eq:estimeigenvDirac}.
Next we prove that the equality in \eqref{eq:estimeigenvDirac} cannot be realized. Assume it were the case, then we would have equalities in all the above inequalities and from Theorem \ref{thm:bessel} the manifold $M$ must be isometric to a geodesic ball.
Furthermore, the spinor field $\psi$ is a Killing spinor (as a consequence of being a twistor spinor and an eigenspinor) with Killing constant $-\frac{\lambda}{n}.$ 
But since on the one hand the scalar curvature of a manifold with such a Killing spinor must be equal to $\frac{4}{n}(n-1)\lambda^2$ and $M$ is Ricci- and hence scalar-flat on the other hand, we deduce that $\lambda=0$.
This contradicts $\mu>0$ and concludes the proof of Theorem \ref{t:estimeigenvDirac}.
\hfill$\square$ 

\noindent Let us now discuss the Dirac spectrum under the MIT bag boundary condition. As we mentioned before, the eigenvalues of the Dirac operator are in this case complex numbers with positive imaginary part.
This fact follows directly from the relation \cite[p. 386]{HMR} 
\begin{equation}\label{identitymit}
2 {\rm Im}(\lambda)\int_M|\psi|^2d\mu_g=\int_{\partial M}|\psi|^2d\mu_g,
\end{equation}
which holds for any eigenspinor $\psi$ associated with an eigenvalue $\lambda$. Now, if we come back to Equality \eqref{Diracschro} with $f=\frac{1}{2}|\psi|^2$, we get after using the nonnegativity of $|P\psi|^2$ that
\begin{eqnarray} \label{laplacianmit}
\Delta f&=&-\frac{S}{2}f+2{\rm Re}(\lambda^2)f-2\frac{|\lambda|^2}{n}f-|P\psi|^2\nonumber\\
&\leq&-\frac{S}{2}f+2({\rm Re}(\lambda^2)-|\lambda|^2)f+\frac{2(n-1)}{n}|\lambda|^2f\nonumber\\
&\leq &\frac{2(n-1)}{n}\left(|\lambda|^2-\frac{n}{4(n-1)}\min_M(S)-\frac{2n}{n-1}{\rm Im}(\lambda)^2\right)f.
\end{eqnarray}
However, we do not have any control on the sign of the r.h.s. of Inequality \eqref{laplacianmit}  in order to deduce an estimate for the eigenvalues using Theorem \ref{thm:bessel} as we did for the other boundary conditions. Notice that S.~Raulot established in \cite[Thm. 1]{R} a lower bound for the first eigenvalue of the Dirac operator with MIT bag condition that involves the imaginary part of $\lambda$ and a lower bound of the mean curvature (assumed to be positive) but unfortunately it still does not provide any new information on the sign of the r.h.s. of \eqref{laplacianmit}.
We can however give a new and short proof of Raulot's estimate \cite[Thm. 1]{R}: 
\begin{thm}[S.~Raulot \cite{R}]\label{thm:raulotestimate} 
Let $(M^n,g)$ be a compact Riemannian spin manifold whose boundary satisfies $H>0.$ Then any eigenvalue $\lambda$ of the Dirac operator of $M$, under the MIT bag boundary condition, satisfies
$$|\lambda|^2\geq \frac{n}{4(n-1)}\min_M(S)+nH_0{\rm Im}(\lambda),$$ 
where $H_0$ is the infimum of the mean curvature. Equality holds if and only if the associated eigenspinor is an imaginary Killing spinor and the boundary is totally umbilical with constant mean curvature.
\end{thm} 

{\noindent \it Proof.} We follow the same steps as we did in Theorem \ref{t:estimeigenvDirac} to get Inequality  \eqref{eq:ineqDeltafDirac}. Indeed, taking into account the boundary condition $i\nu\cdot\psi=\psi$ on $\partial M,$ we compute
\begin{eqnarray*}
\int_M\Delta f d\mu_g&\bui{=}{\rm\eqref{eq:DiracMdM}}&\int_{\partial M}\Re(\langle-\nu\cdot D\psi-D^{\partial M}\psi+\frac{(n-1)H}{2}\psi,\psi\rangle)d\mu_g\\
&=&\int_{\partial M}\Re(-\lambda\langle\nu\cdot\psi,\psi\rangle)d\mu_g-\int_{\partial M}\langle D^{\partial M}\psi,\psi\rangle d\mu_g+\frac{n-1}{2}\int_{\partial M}H|\psi|^2d\mu_g\\
&\geq &((n-1)H_0-2 {\rm Im}(\lambda))\int_{\partial M}f d\mu_g.
\end{eqnarray*}
Here, as before we use the fact that $\displaystyle\int_{\partial M}\langle D^{\partial M}\psi,\psi\rangle d\mu_g\leq0$ which is also valid for the MIT bag boundary condition \cite{HMR}. 
Integrating Inequality \eqref{laplacianmit} over $M$ yields the desired inequality after using the identity \eqref{identitymit}. If now equality holds, then the eigenspinor is a Killing spinor of Killing number $\frac{-\lambda}{n}.$ But as $\lambda$ is a complex number (remember that its imaginary part is positive), it must be purely imaginary (see e.g. \cite[Ch. 8]{BHMM} or \cite[Sec. A.4]{GinsurveyDirac}) which implies that $\psi$ is an imaginary Killing spinor. The last part follows from differentiating the boundary condition along any vector field tangent to the boundary, see e.g. \cite[pp. 142-143]{R}.
This finishes the proof.
\hfill$\square$

\noindent As we see from Theorem \ref{thm:raulotestimate}, if one assumes $H_0>\frac{2}{n-1}{\rm Im}(\lambda),$ then the r.h.s. of Inequality \eqref{laplacianmit} is in this case positive and therefore Theorem \ref{thm:bessel} can be applied. However, we think that it is unnatural to require such a bound on the mean curvature as it depends on the eigenvalue $\lambda$ in question. 

\noindent Another way for estimating the eigenvalues of the Dirac operator is to look at a conformal class of metrics, we refer to e.g. \cite[Sec. 2.3 \& 5.4]{BHMM} or \cite[Sec. 3.3]{GinsurveyDirac} for general facts on the subject.
In this case, the spectrum of the Dirac operator is known to be related to the spectrum of the so-called {\it Yamabe operator} through the so-called Hijazi estimate \cite{Hij86}.
In \cite{R1} S.~Raulot proved  that, under the CHI or the MIT bag condition, any eigenvalue $\lambda$ of the Dirac operator satisfies, for $n\geq 3,$ 
\begin{equation}\label{hijazi}
|\lambda|^2\geq \frac{n}{4(n-1)}\mu_1(Y),
\end{equation}
where the inequality is strict for the MIT bag condition and  characterizes in its limiting case the half round sphere for the CHI condition. Here, $\mu_1(Y)$ denotes the first positive eigenvalue of the Yamabe problem originally defined by Escobar in  \cite{E}:
\begin{equation}\label{yamabe}
\left\{
\begin{array}{lll}
\medskip
Y(f):=\frac{4(n-1)}{n-2}\Delta f+S f=\mu_1(Y) f\, && \text{on }\,\, M, \\
\medskip
\frac{\partial f}{\partial\nu}=\frac{n-2}{2}Hf && \text{on}\,\, \partial M.\\
\end{array}
\right.
\end{equation} 
Recall that $\nu$ denotes the inward unit normal vector field along $\partial M$.
If the mean curvature is nonnegative then it is easy to check, after multiplying the first equation in \eqref{yamabe} involving a first eigenfunction $f$ with $f$ itself and integrating over $M$, that the inequality $\mu_1(Y)\geq \mathop{{\rm min}}\limits_{M}(S)$ holds, with equality if and only if $S$ is constant on $M$ and $H=0$ on $\partial M$.
We will use this last fact to apply Theorem \ref{thm:bessel} to an eigenfunction of the Yamabe operator in order to deduce an estimate for $\mu_1(Y)$ in terms of zeros of Bessel functions.
This will allow later to derive a new estimate for the Dirac operator.

\begin{thm}\label{thm:yamabe}
Let $(M^n,g)$ be a compact Riemannian manifold of dimension $n\geq 3$ with smooth boundary.
Assume that the Ricci curvature of $M$ is nonnegative and the mean curvature is bounded from below by $H_0>0.$
Let $\tau_1$ be the only positive zero of $x\mapsto x\frac{J_{\frac{n}{2}(x)}}{J_{\frac{n}{2}-1}(x)}-\frac{n-2}{2}$ on $]0,j_{\frac{n}{2}-1,1}[$. Then 
\begin{equation}\label{eq:minormu1Yamabe}\mu_1(Y)\geq \min_M(S)+\frac{4(n-1)}{n-2}\tau_1^2H_0^2.\end{equation}
Equality is realized if and only if the manifold $M$ is isometric to a round ball in $\R^n$.
\end{thm}

{\noindent \it Proof.} Let $f$ be an eigenfunction of the problem \eqref{yamabe} associated with the eigenvalue $\mu_1(Y)$.
Recall that $f$ cannot change its sign, so that $f$ can be assumed to be positive in the interior of $M$.
Then we have
$$\Delta f=\frac{n-2}{4(n-1)}\left(\mu_1(Y)-S\right)f\leq \frac{n-2}{4(n-1)}\left(\mu_1(Y)-\min_M(S)\right)f.$$ 
Let $\mu:=\frac{n-2}{4(n-1)}\left(\mu_1(Y)-\min_M(S)\right).$
Notice that $\mu>0$ since $H$ cannot vanish.
If by contradiction $\frac{\sqrt{\mu}}{H_0}<\tau_1,$ then from Theorem \ref{thm:bessel} we deduce that 
\begin{equation}\label{eq:intfYamabe}\int_{\partial M}fd\mu_g\geq\sqrt{\mu}\frac{J_{\frac{n}{2}-1}(\frac{\sqrt{\mu}}{H_0})}{J_{\frac{n}{2}}(\frac{\sqrt{\mu}}{H_0})}\int_M fd\mu_g.\end{equation}
Note that this implies that $\int_{\partial M}f d\mu_g\neq0$.
Integrating the inequality $\Delta f\leq \mu f$ along with 
$$\int_M\Delta f d\mu_g=\int_M\frac{\partial f}{\partial\nu} d\mu_g\geq \frac{n-2}{2}H_0\int_{\partial M} f d\mu_g,$$
yields $\frac{2}{n-2}\geq\frac{J_{\frac{n}{2}-1}(\frac{\sqrt{\mu}}{H_0})}{\frac{\sqrt{\mu}}{H_0}J_{\frac{n}{2}}(\frac{\sqrt{\mu}}{H_0})}.$ Finally, as we did in the proof of Theorem \ref{t:estimeigenvDirac}, we use the fact that the function $x\mapsto \frac{J_{\frac{n}{2}-1}(x)}{xJ_{\frac{n}{2}}(x)}$ is decreasing on $]0,j_{\frac{n}{2}-1,1}[$ to get the contradiction. 
This proves \eqref{eq:minormu1Yamabe}.
If \eqref{eq:minormu1Yamabe} is an equality, then because of $\frac{\sqrt{\mu}}{H_0}=\tau_1\in\,]0,j_{\frac{n}{2}-1,1}[$, the inequality \eqref{eq:intfYamabe} still applies and must be an equality, therefore $M$ must be isometric to a round ball in $\R^n$ by Theorem \ref{thm:bessel}.
Conversely, if $M$ is a round ball in $\R^n$ of radius $\frac{1}{H_0}$, then the two inequalities involving $\Delta f$ and used in the proof of \eqref{eq:minormu1Yamabe} are equalities since scalar and mean curvatures are constant.
Moreover, for a round ball in $\R^n$, the problem \eqref{yamabe} reduces to the Robin boundary value problem \eqref{eq:robinproblem} for the first eigenvalue $\mu:=\frac{n-2}{4(n-1)}\mu_1(Y)$ and where $\tau:=\frac{n-2}{2}H_0>0$.
As was already noticed in Remark \ref{rem:eigenvalueball}, the first eigenvalue of the Robin boundary value problem on a round ball of $\R^n$ always satisfies $\frac{\sqrt{\mu}}{H_0}\in\,]0,j_{\frac{n}{2}-1,1}[\subset\,]0,j_{\frac{n}{2},1}[$, therefore \eqref{eq:intfYamabe} applies and is actually an equality again by Theorem \ref{thm:bessel}.
On the whole, all three inequalities used in the proof of \eqref{eq:minormu1Yamabe} are equalities for a round ball, therefore \eqref{eq:minormu1Yamabe} itself must be an equality.
This shows the equivalence for the limiting case and concludes the proof of Theorem \ref{thm:yamabe}.
\hfill$\square$

\noindent Combining Inequality \eqref{hijazi} with \eqref{eq:minormu1Yamabe}, we deduce the following
\begin{cor}\label{cor:chiestimate}
Let $(M^n,g)$ be a compact Riemannian spin manifold of dimension $n\geq 3$ with smooth boundary.
Assume that the Ricci curvature of $(M^n,g)$ is nonnegative and the mean curvature of $\partial M$ is bounded from below by $H_0>0.$
Let $\tau_1$ be the only positive zero of $x\mapsto x\frac{J_{\frac{n}{2}(x)}}{J_{\frac{n}{2}-1}(x)}-\frac{n-2}{2}$ on $]0,j_{\frac{n}{2}-1,1}[$. Then, under the CHI or the MIT bag conditions, any eigenvalue $\lambda$ of the Dirac operator satisfies  
\begin{equation}\label{eq:minorDiracmu1}|\lambda|^2>\frac{n}{4(n-1)}\min_M(S)+\frac{n}{n-2}\tau_1^2H_0^2.\end{equation}
\end{cor}
Note that equality cannot hold in \eqref{eq:minorDiracmu1} since  for the MIT bag boundary condition \eqref{hijazi} is anyway strict while for the CHI boundary condition equality in \eqref{hijazi} implies minimality of $\partial M$ in $M$.

\noindent As we can see from Theorem \ref{t:estimeigenvDirac} and Corollary \ref{cor:chiestimate} that there are two different but analogous estimates for the first eigenvalue of the Dirac operator under the CHI boundary condition. 
One might ask if there is a way to compare the numbers $\frac{n}{2(n-1)}\tau_0^2$ and $\frac{n}{n-2}\tau_1^2$ in order to check which estimate is better. Recall here that $\tau_0$ and $\tau_1$ are respectively the first positive zeros of the functions $x\mapsto x\frac{J_{\frac{n}{2}(x)}}{J_{\frac{n}{2}-1}(x)}-(n-1)$ and $x\mapsto x\frac{J_{\frac{n}{2}(x)}}{J_{\frac{n}{2}-1}(x)}-\frac{n-2}{2}$ on $]0,j_{\frac{n}{2}-1,1}[$.
It turns out that \cite{Freitasprivate}, for any $n\geq3$,
\begin{equation}\label{eq:ineqFreitastau0tau1}\left(\frac{\tau_1}{\tau_0}\right)^2\geq\frac{(n+1)(n+1-\sqrt{4n+1})}{n(n-1)},\end{equation}
which implies that $\left(\frac{\tau_1}{\tau_0}\right)^2>\frac{n-2}{2(n-1)}$, i.e. $\frac{n}{n-2}\tau_1^2>\frac{n}{2(n-1)}\tau_0^2$.
Therefore, \eqref{eq:minorDiracmu1} is better than \eqref{eq:estimeigenvDirac} for any $n\geq3$.
Moreover, because of $\tau_1<\tau_0$, inequality \eqref{eq:ineqFreitastau0tau1} implies that $\frac{\tau_1}{\tau_0}\underset{n\to\infty}{\longrightarrow}1$, which means that $\frac{\frac{n}{n-2}\tau_1^2}{\frac{n}{2(n-1)}\tau_0^2}\underset{n\to\infty}{\longrightarrow}2$: asymptotically as $n\to\infty$, the lower bound in \eqref{eq:minorDiracmu1} is even much better than \eqref{eq:estimeigenvDirac}.
This could be explained by the fact that less information is lost upon proving \eqref{hijazi} and \eqref{eq:minormu1Yamabe} than proving \eqref{eq:estimeigenvDirac} directly.


\section{Robin Laplacian for differential forms}\label{s:RobinLapldiffforms}


We first recall the so-called Lopatinski\u{\i}-Shapiro criterion for ellipticity of boundary value problems, see e.g. \cite[Sec. 1.6]{SchwarzLNM95} to which we shall stay close.
Let $(M^n,g)$ be any Riemannian manifold with nonempty boundary $\partial M$.
Let $P$ be any $k^{\textrm{th}}$-order linear differential operator acting on sections of some Riemannian or Hermitian vector bundle $E\to M$.
A \emph{boundary condition} will be considered here as the direct sum $\displaystyle{\bigoplus_{j=1}^l B_j}$ of linear differential operators $B_j\colon\Gamma(M,E)\to\Gamma(\partial M,E_j)$ of order $k_j<k$, where $E_j\to\partial M$, $1\leq j\leq l$, are Riemannian or Hermitian vector bundles.
We consider the following boundary value problem: for any $f\in\Gamma(M,E)$ and $u_j\in\Gamma(\partial M,E_j)$, $1\leq j\leq l$, find $u\in\Gamma(M,E)$ solving
\begin{equation}\label{eq:BVPgeneral}
\left\{\begin{array}{lll}Pu&=f&\textrm{ on }M\\ B_ju&=u_j&\textrm{ on }\partial M,\;\forall\,1\leq j\leq l\end{array}\right.
\end{equation}
Let $\sigma_P$ and $\sigma_{B_j}$ be the principal symbols of the operators $P$ and $B_j$ respectively, $1\leq j\leq l$.
In our convention, for any smooth function $f$ defined in a neighbourhood of a point $x$,
\[\sigma_P(d_xf):=\underbrace{\left[\cdots[[P,f],f],\cdots,f\right]}_{k\textrm{ pairs of brackets}}\in\mathrm{Hom}(E_x,E_x)\]
and analogously for $B_j$; of course, $\sigma_B=B$ if $B$ has vanishing order.
In order to formulate the Lopatinski\u{\i}-Shapiro ellipticity condition, the following space must be defined: given any $x\in\partial M$ and $v\in T_x\partial M$, let
\[\mathcal{M}_{v}^+:=\left\{\textrm{bounded solutions }y=y(t)\textrm{ on }\R_+\textrm{ to the ODE }\sigma_P((-iv,\partial_t))y=0\right\}.\]
Here the map $\sigma_P((-iv,\partial_t))$ must be understood as follows: considering $\sigma_P$ pointwise as a homogeneous polynomial of degree $k$ on $T_xM$, we apply it to the one-form $-iv^\flat+\partial_t\cdot\nu^\flat$, where $\nu$ is the inner unit normal at $x$ and where we see $\partial_t$ as a coefficient; what we obtain at the end is a $k^{\textrm{th}}$-order linear differential operator in one variable and with constant coefficients.
In particular the ODE $\sigma_P((-iv,\partial_t))y=0$ has a $k$-dimensional space of solutions that are defined on $\R$. The following definition is taken from \cite[Def. 1.6.1]{SchwarzLNM95}:
\begin{definition}\label{def:LopatinskiiShapiroellipticbc}
The boundary value problem {\rm(\ref{eq:BVPgeneral})} is called \emph{elliptic} if and only if both following conditions are satisfied:
\begin{enumerate}[(a)]
\item The differential operator $P$ is itself elliptic, that is, for any $x\in M$ and $\xi\in T_x^*M$, the map $\sigma_P(\xi)\colon E_x\to E_x$ is an isomorphism.
\item For all $x\in \partial M$ and $v\in T_x\partial M$, the map 
\begin{eqnarray*}
\mathcal{M}_{v}^+&\longrightarrow&\bigoplus_{j=1}^l (E_j)_x\\
y&\longmapsto&\left(\sigma_{B_1}((-iv,\partial_t))y,\ldots,\sigma_{B_l}((-iv,\partial_t))y\right)(0)
\end{eqnarray*}
is an isomorphism.
\end{enumerate}
\end{definition}
Now we look at the following setting where $E=\Lambda^pT^*M$, $E_1=\Lambda^pT^*\partial M$, $E_2=\Lambda^{p-1}T^*\partial M$, $P=\Delta=d\delta_g+\delta_g d$, $B_1\omega=\iota^*(\nu\lrcorner \,d\omega-\tau\omega)$ and $B_2\omega=\iota^*(\nu\lrcorner\,\omega)$ for any given $p\in\{0,1,\ldots,n\}$, $\tau\in\R$ and any $\omega\in\Omega^p(M)$.
Recall that $\iota\colon\partial M\to M$ denotes the inclusion map.
The principal symbols are given by 
\begin{eqnarray*}
\sigma_P(\xi)\omega&=&-|\xi|^2\cdot\omega\\
\sigma_{B_1}(\xi)\omega&=&\xi_\nu\cdot \iota^*\omega-\xi_T\wedge \iota^*(\nu\lrcorner\,\omega)\\
\sigma_{B_2}(\xi)\omega&=&\iota^*(\nu\lrcorner\,\omega),
\end{eqnarray*}
where $\xi=\xi_\nu\nu^\flat+\xi_T\in\R\cdot\nu^\flat\oplus T_x^*\partial M$ for any $\xi\in T_x^*M$ and $x\in\partial M$.

\begin{thm}\label{p:Robinisellipticandselfadjoint} Let $(M,g)$ be a compact Riemannian manifold with smooth boundary.
Fix a positive real number $\tau$ and consider the eigenvalue problem
\begin{equation}\label{eq:BVPRobin}
\left\{\begin{array}{lll}\Delta \omega&=\lambda\omega &\textrm{ on }M\\ \iota^*(\nu\lrcorner\,d\omega-\tau\omega)&=0&\textrm{ on }\partial M\\
\iota^*(\nu\lrcorner\,\omega)&=0&\textrm{ on }\partial M.\end{array}\right.
\end{equation}
Then, we have 
\begin{enumerate}
\item\label{claim:Robinellselfadj} The boundary value problem {\rm(\ref{eq:BVPRobin})} is elliptic in the sense of {\rm Definition \ref{def:LopatinskiiShapiroellipticbc}} and self-adjoint. As a consequence, it admits an increasing unbounded sequence of nonnegative real eigenvalues with finite multiplicities $\lambda_{1,p}(\tau)<\lambda_{2,p}(\tau)\leq\cdots$
\item\label{claim:Robinpostaupos}{Actually $\lambda_{1,p}(\tau)>0$ holds, i.e. \eqref{eq:BVPRobin} has trivial kernel.}
\end{enumerate}
\end{thm}

{\noindent \it Proof.}
The Laplace operator on forms is clearly elliptic since for any nonvanishing $\xi\in T^*M$ the map $-|\xi|^2\cdot\mathrm{Id}$ is an isomorphism.
Moreover, for any $v\in T^*\partial M$, we can write $\sigma_{\Delta}((-iv,\partial_t))=-\langle(-iv,\partial_t),(-iv,\partial_t)\rangle=|v|^2-\partial_t^2$, so that
\[\mathcal{M}_{v}^+=\left\{e^{-t|v|}\cdot\omega_0,\;\omega_0\in\Lambda^pT_x^*M\right\}.\]
On the other hand, $\sigma_{B_1}((-iv,\partial_t))\omega=\partial_t(\iota^*\omega)+iv^\flat\wedge  \iota^*(\nu\lrcorner\,\omega)$ and $\sigma_{B_2}((-iv,\partial_t))\omega=\iota^*(\nu\lrcorner\,\omega)$, so that if $y=e^{-t|v|}\cdot\omega_0$ is any element of $\mathcal{M}_{v}^+$, then 
\[\left(\sigma_{B_1}((-iv,\partial_t))y,\sigma_{B_2}((-iv,\partial_t))y\right)(0)=\left(-|v|\iota^*\omega_0+iv^\flat\wedge  \iota^*(\nu\lrcorner\,\omega_0),\iota^*(\nu\lrcorner\,\omega_0)\right).\]
If the r.h.s. of that identity vanishes, then $\iota^*(\nu\lrcorner\,\omega_0)=0$ and therefore $|v|\iota^*\omega_0=0$, which under the assumption $v\neq0$ yields $\iota^*\omega_0=0$ and thus $\omega_0=0$.
This shows that the map $\mathcal{M}_{v}^+\longrightarrow\bigoplus_{j=1}^2(E_j)_x$ of Definition \ref{def:LopatinskiiShapiroellipticbc} is injective and hence an isomorphism by equality of the space dimensions.
This shows (\ref{eq:BVPRobin}) to be elliptic in the sense of Definition \ref{def:LopatinskiiShapiroellipticbc}. To show self-adjointness, we need to prove that, for any compactly-supported smooth $p$-forms $\omega,\omega'$ on $M$ satisfying the boundary conditions $B_1\omega=B_2\omega=B_1\omega'=B_2\omega'=0$, the identity $\displaystyle{\int_M\langle\Delta\omega,\omega'\rangle d\mu_g=\int_M\langle\omega,\Delta\omega'\rangle d\mu_g}$ holds.
But using the partial integration formula
\[\int_M\langle d\alpha,\beta\rangle d\mu_g=\int_M\langle\alpha,\delta_g\beta\rangle d\mu_g-\int_{\partial M}\langle \iota^*\alpha,\nu\lrcorner\,\beta\rangle d\mu_g\]
that is valid for all $\alpha\in\Omega^p(M)$ and $\beta\in\Omega^{p+1}(M)$, we obtain
\begin{eqnarray}\label{eq:laplacian}
\int_M\langle\Delta\omega,\omega'\rangle d\mu_g&=&\int_M\langle d\delta_g\omega,\omega'\rangle d\mu_g+\int_M\langle \delta_gd\omega,\omega'\rangle d\mu_g\nonumber\\
&=&\int_M\langle\delta_g\omega,\delta_g\omega'\rangle d\mu_g-\int_{\partial M}\langle \iota^*\delta_g\omega,\nu\lrcorner\,\omega'\rangle d\mu_g\nonumber\\
&&+\int_M\langle d\omega,d\omega'\rangle d\mu_g+\int_{\partial M}\langle\nu\lrcorner d\omega,\iota^*\omega'\rangle d\mu_g\nonumber\\
&=&\int_M\langle d\omega,d\omega'\rangle d\mu_g+\langle\delta_g\omega,\delta_g\omega'\rangle d\mu_g\nonumber\\&&+\int_{\partial M}\langle \iota^*(\nu\lrcorner d\omega),\iota^*\omega'\rangle d\mu_g-\langle \iota^*\delta_g\omega,\underbrace{\iota^*(\nu\lrcorner\,\omega')}_{0}\rangle  d\mu_g\\
&=&\int_M\langle d\omega,d\omega'\rangle d\mu_g+\langle\delta_g\omega,\delta_g\omega'\rangle d\mu_g+\int_{\partial M}\tau \langle \iota^*\omega,\iota^*\omega'\rangle d\mu_g,\nonumber
\end{eqnarray}
which is clearly symmetric in $(\omega,\omega')$ because of $\tau\in\R$. This shows (\ref{eq:BVPRobin}) to be self-adjoint. As a consequence, the spectrum of the Robin operator consists of an unbounded sequence of real eigenvalues of finite multiplicities.
Moreover, if $\tau>0$, then for any differential $p$-form $\omega$ on $M$, we can deduce from the above computation the equality 
\begin{equation}\label{eq:laplacianomega2}
\int_M \langle \Delta\omega,\omega\rangle d\mu_g= \int_M ( |d\omega|^2+|\delta_g\omega|^2) d\mu_g+\tau\int_{\partial M} |\iota^*\omega|^2d\mu_g \geq 0.
\end{equation}
Therefore, the spectrum of the Robin operator for $\tau>0$ must be nonnegative and can therefore be written as an increasing unbounded sequence $\lambda_{1,p}(\tau)\leq\lambda_{2,p}(\tau)\leq\ldots$ of nonnegative real eigenvalues of finite multiplicities.
Note that, by Courant's nodal domain theorem, the first eigenvalue $\lambda_{1,p}(\tau)$ is simple and every associated eigenfunction cannot change its sign on $M$.
This shows claim \ref{claim:Robinellselfadj}. Next we show that $0$ is not an eigenvalue when $\tau>0$. Let $\omega$ lie in the kernel of $\Delta$. From the formula above, we obtain $d\omega=\delta_g\omega=0$ on $M$ and $\iota^*\omega=0$ on $\partial M.$
But using the identity $|\omega|^2=|\iota^*\omega|^2+|\nu\lrcorner \omega|^2$ at any point on the boundary, we deduce that $\omega=0$ on $\partial M.$ 
{ Now by \cite[Thm. p. 445]{Anne89}, any harmonic form on $M$ that vanishes along $\partial M$ must vanish identically, therefore $\omega=0$.
This proves claim \ref{claim:Robinpostaupos} and concludes the proof.\hfill$\square$
}

\begin{prop}\label{claim:varcharfirsteigenvRobin} Let $(M,g)$ be a compact Riemannian manifold with smooth boundary. Let $\tau$ be a positive number.
Then the first eigenvalue $\lambda_{1,p}(\tau)$ of the Robin boundary problem {\rm(\ref{eq:BVPRobin})} can be characterized as follows:  
$$\lambda_{1,p}(\tau)={\rm inf}\left\{\frac{\displaystyle\int_M \left(|d\omega|^2 +|\delta_g\omega|^2\right)d\mu_g+\tau\int_{\partial M} |\iota^*\omega|^2d\mu_g}{\displaystyle\int_M|\omega|^2d\mu_g}\right\},$$
where $\omega$ runs over all non-identically vanishing $p$-forms on $M$ such that $\nu\lrcorner\omega=0$.
\end{prop}

{\noindent \it Proof.} We have seen that Identity \eqref{eq:laplacianomega2} holds as soon as $\omega$ satisfies the boundary condtions $\iota^*(\nu\lrcorner d\omega)=\tau\iota^*\omega$ and $\iota^*(\nu\lrcorner\omega)=0$. In particular, this proves the variational characterization 
\[\lambda_{1,p}(\tau)=\inf_{\substack{\omega\in\Omega^p(M)\setminus\{0\}\\ \iota^*(\nu\lrcorner d\omega)=\tau \iota^*\omega,\;\iota^*(\nu\lrcorner\omega)=0}}\left\{\frac{\int_M(|d\omega|^2+|\delta_g\omega|^2)d\mu_g+\tau\int_{\partial M}|\iota^*\omega|^2 d\mu_g}{\int_M|\omega|^2d\mu_g}\right\}.\]
We next show that the boundary condition $\iota^*(\nu\lrcorner d\omega)=\tau \iota^*\omega$ can actually be dropped off in the infimum above.
We follow \cite[Ch. 5, Sec. 9]{TaylorPDE1}.
Define the standard Sobolev spaces:
$$H^k(M,\Lambda^p):=\{\omega\in L^2(M,\Lambda^pT^*M)|\,\nabla^l\omega\in L^2(M,\otimes^lT^*M\otimes\Lambda^pT^*M), \,\,\textrm{for any}\,\, 1\leq l\leq k\}.$$
As in \cite[Sec. 5.9]{TaylorPDE1}, we define the following closed subspaces of $H^k(M,\Lambda^p)$:
\begin{eqnarray*} 
H_{\rm Rob}^1(M,\Lambda^p)&:=&\{\omega\in H^1(M,\Lambda^p)|\,\iota^*(\nu\lrcorner\omega)=0\textrm{ on }\partial M\}\\
H_{\rm Rob}^2(M,\Lambda^p)&:=&\{\omega\in H^2(M,\Lambda^p)|\,\iota^*(\nu\lrcorner\omega)=0\textrm{ and }\iota^*(\nu\lrcorner d\omega)-\tau \iota^*\omega=0\textrm{ on }\partial M\}.
\end{eqnarray*}
Note that the conditions defining those subspaces make sense because of the existence of a continuous extension $H^1(M,\Lambda^p)\to L^2(\partial M,\Lambda^p)$ of the trace map $\omega\mapsto\omega_{|_{\partial M}}$.
Consider, for any $\omega\in H_{\rm Rob}^1(M,\Lambda^p)$, the linear operator 
\begin{equation}\label{eqlrobin}
\omega'\longmapsto (L_{\rm Rob}\omega)(\omega'):=\int_M\left(\langle d\omega,d\omega'\rangle+\langle\delta_g\omega,\delta_g\omega'\rangle\right) d\mu_g+\tau\int_{\partial M}\langle \iota^*\omega,\iota^*\omega'\rangle d\mu_g
\end{equation}
on $H_{\rm Rob}^1(M,\Lambda^p)$. Note that this operator just differs from the one in \cite[Eq. 9.21]{TaylorPDE1} by the boundary term involving $\tau$. 
 Clearly, $L_{\rm Rob}$ defines a bounded linear operator 
\[L_{\rm Rob}\colon H_{\rm Rob}^1(M,\Lambda^p)\longrightarrow H_{\rm Rob}^1(M,\Lambda^p)^*,\]
where $H_{\rm Rob}^1(M,\Lambda^p)^*$ denotes the topological dual of $H_{\rm Rob}^1(M,\Lambda^p)$. This can be proved by estimating the boundary term in \eqref{eqlrobin} using Cauchy-Schwarz inequality and the trace theorem. Moreover, there exists a positive constant $C_0$ such that, for every $\omega\in H_{\rm Rob}^1(M,\Lambda^p)$,
\begin{equation}\label{eq:LRobbddbelow}
\left((L_{\rm Rob}+C_0\right)\omega,\omega)\geq C\|\omega\|_{H^1}^2
\end{equation}
for some further positive constant $C.$ Inequality \eqref{eq:LRobbddbelow} follows in a straightforward way from \cite[Eq. 9.24]{TaylorPDE1} since the boundary term coming from \eqref{eqlrobin} is  positive because of $\tau>0$ by assumption. Now Inequality \eqref{eq:LRobbddbelow} not only shows that $L_{\rm Rob}+C_0$ is injective with closed range but also that it is bijective, see e.g. \cite[Prop. 9.5]{TaylorPDE1}. We denote by $T_{\rm Rob}:H_{\rm Rob}^1(M,\Lambda^p)^*\to H_{\rm Rob}^1(M,\Lambda^p)$ the inverse map as well as the induced map $T_{\rm Rob}:L^2(M,\Lambda^p)\to L^2(M,\Lambda^p)$. Note here that there is a compact embedding $H_{\rm Rob}^1(M,\Lambda^p)\to L^2(M,\Lambda^p),$ since the embedding $H^1(M,\Lambda^p)\to L^2(M,\Lambda^p)$ is already compact. The operator $T_{\rm Rob}$ on the $L^2$-level is compact and, being selfadjoint and positive, has a discrete spectrum which can be described as a nonincreasing sequence of positive eigenvalues of finite multiplicity converging to $0$. Moreover, there exists an $L^2$-orthonormal basis $(\omega_j)_{j\in\mathbb{N}}$ of $L^2(M,\Lambda^p)$ consisting of eigenvectors for $T_{\rm Rob}$: for every $j\in\mathbb{N}$, we have $T_{\rm Rob}\omega_j=\mu_j\omega_j$.
Note that necessarily $\omega_j\in H_{\rm Rob}^1(M,\Lambda^p)$ holds for every $j$ since the range of $T_{\rm Rob}$ actually lies in $H^1_{\rm Rob}(M,\Lambda^p)$ by definition.
By construction, $(\omega_j)_{j\in\mathbb{N}}$ is an $L^2$-orthonormal basis of eigenvectors for $L_{\rm Rob}$ since $L_{\rm Rob}\omega_j=(\frac{1}{\mu_j}-C_0)\omega_j$ holds for every $j$. The central remark is now the following: for every $\omega\in L^2(M,\Lambda^p)$, actually $u:=T_{\rm Rob}\omega\in H_{\rm Rob}^2(M,\Lambda^p)$ must hold. To see this, we shall divide the proof into several steps whose technical details will be ignored since they are completely analogous to those from the proof of \cite[Prop. 9.6]{TaylorPDE1}.\\
$\bullet$ {\rm Step} 1: Due to the ellipticity of the operator $L_{\rm Rob}$ (see e.g. \cite[Prop. 7.2]{TaylorPDE1} for the corresponding estimate), we have that $u=T_{\rm Rob}\omega$  lies in $H^2(M,\Lambda^p).$\\ 
$\bullet$ {\rm Step} 2: We must prove that the boundary condition $\iota^*(\nu\lrcorner du)-\tau \iota^*u=0$ holds for $u$ if and only if the boundary term 
$$\int_{\partial M}\langle \iota^*(\nu\lrcorner du)-\tau \iota^*u,\iota^*\omega'\rangle  d\mu_g$$ 
vanishes for all $\omega'\in H_{\rm Rob}^1(M,\Lambda^p)$. Note that this is not obvious since {\sl a priori} the range of the trace map restricted to $H_{\rm Rob}^1(M,\Lambda^p)$ is not a dense subspace of $L^2(\partial M,\Lambda^p)$. To prove the if condition, we extend the pointwise homomorphism field $\nu\lrcorner\;\colon\Lambda^{p+1}T^*M_{|_{\partial M}}\to\Lambda^pT^*M_{|_{\partial M}}$ along $\partial M$ to a smooth homomorphism field $\sigma\colon \Lambda^{p+1}T^*M\to\Lambda^pT^*M$ on $M$.
Pick any $\alpha\in\Omega^{p+1}(M)$ and put $\omega':=\sigma\alpha\in\Omega^p(M)$.
Note that $\iota^*(\nu\lrcorner\omega')=\iota^*(\nu\lrcorner\nu\lrcorner\alpha)=0$ holds along $\partial M$, therefore $\omega'\in H_{\rm Rob}^1(M,\Lambda^p)$.
Moreover,
\[\int_{\partial M}\langle \iota^*(\nu\lrcorner du)-\tau \iota^*u,\iota^*\omega'\rangle  d\mu_g=0\;\Longleftrightarrow\;\int_{\partial M}\langle \nu\wedge\left(\iota^*(\nu\lrcorner du)-\tau \iota^*u\right),\iota^*\alpha\rangle  d\mu_g=0.\]
This holds for all $\alpha\in\Omega^{p+1}(M)$, therefore $\nu\wedge\left(\iota^*(\nu\lrcorner du)-\tau \iota^*u\right)=0$ along $\partial M$.
Now taking the interior product with $\nu$ allows to deduce that $\iota^*(\nu\lrcorner du)-\tau \iota^*u=0$, as claimed.\\
$\bullet$ {\rm Step} 3: We want to show that $\Delta(T_{\rm Rob}\omega)=\omega-C_0 T_{\rm Rob}\omega$ holds in the distributional sense. Recall that $T_{\rm Rob}=(L_{\rm Rob}+C_0)^{-1}$. Taking the scalar product of the relation $(L_{\rm Rob}+C_0)(u)=\omega$ with any $\omega'\in H_{\rm Rob}^1(M,\Lambda^p),$  we get 
\begin{eqnarray*}
\int_M \langle \omega,\omega'\rangle d\mu_g&=&\int_M \langle L_{\rm Rob} u,\omega'\rangle d\mu_g+C_0 \int_M\langle   u,\omega'\rangle d\mu_g\\
&=& \int_M (\langle du,d\omega'\rangle+\langle \delta_g u,\delta_g\omega'\rangle) d\mu_g+\tau\int_{\partial M}\langle\iota^*u,\iota^*\omega'\rangle+C_0 \int_M\langle  u,\omega'\rangle d\mu_g\\
&\bui{=}{\eqref{eq:laplacian}}& \int_M \langle \Delta u,\omega'\rangle d\mu_g-\int_{\partial M}\langle \iota^*(\nu\lrcorner du),\iota^*\omega'\rangle d\mu_g+\tau\int_{\partial M}\langle\iota^*u,\iota^*\omega'\rangle d\mu_g\\&&+C_0 \int_M\langle  u,\omega'\rangle d\mu_g.
\end{eqnarray*}
This means that
\begin{equation}\label{eq:equationtrobin}
 \int_M \langle (\Delta+C_0)T_{\rm Rob}\omega,\omega'\rangle d\mu_g=\int_M \langle \omega,\omega'\rangle d\mu_g+\int_{\partial M}\langle \iota^*(\nu\lrcorner du-\tau u),\iota^*\omega'\rangle d\mu_g.
 \end{equation}
But this holds for any $\omega'\in H_{\rm Rob}^1(M,\Lambda^p).$ In particular, if we choose $\omega'\in \Omega^{p+1}(M)$ with support away from $\partial M$ (see \cite[p. 407]{TaylorPDE1}) we deduce that $\Delta(T_{\rm Rob}\omega)=\omega-C_0 T_{\rm Rob}\omega$ as required. This ensures the boundary term in \eqref{eq:equationtrobin} has to vanish for all $\omega'\in H_{\rm Rob}^1(M,\Lambda^p)$. Therefore by Step 2 the boundary condition $\iota^*(\nu\lrcorner du)-\tau \iota^*u=0$ must be fulfilled.\\
To conclude the proof, every eigenvector $\omega_j$ of $L_{\rm Rob}$ is an eigenvector for the Laplace operator on $M$ and, since it belongs to the range of $T_{\rm Rob}$, it must lie in $H_{\rm Rob}^2(M,\Lambda^p)$, in particular satisfy the boundary condition of first order.

\hfill$\square$

\noindent From now, unless otherwise stated, we assume $\tau>0$. Recall the absolute boundary conditions 
\begin{equation}\label{eq:absolutebc}
\left\{\begin{array}{lll}\Delta \omega&=\lambda^N\omega &\textrm{ on }M\\ \iota^*(\nu\lrcorner\,d\omega)&=0&\textrm{ on }\partial M\\
\iota^*(\nu\lrcorner\,\omega)&=0&\textrm{ on }\partial M\end{array}\right.
\end{equation}
which generalize the Neumann boundary problem for functions. The spectrum of this Laplacian is discrete and consists of eigenvalues $(\lambda^N_{i,p})_i$ such that $\lambda_{1,p}^N\leq\lambda_{2,p}^N\leq\cdots.$ The Hodge star operator exchanges the absolute boundary conditions and the relative ones which are given by 
\begin{equation}\label{eq:relativebc}
\left\{\begin{array}{lll}\Delta \omega&=\lambda^{R}\omega &\textrm{ on }M\\ \iota^*\omega &=0&\textrm{ on }\partial M\\
\iota^*(\delta_g\omega)&=0&\textrm{ on }\partial M.\end{array}\right.
\end{equation}

\noindent By the min-max principle, the first eigenvalue $\lambda^N_{1,p}$ of the Laplacian is characterized by 
$$\lambda^N_{1,p}={\rm inf}\left\{\frac{\displaystyle\int_M (|d\omega|^2+|\delta_g\omega|^2)d\mu_g}{\displaystyle\int_M|\omega|^2 d\mu_g}\right\},$$
where $\omega$ runs over all $p$-forms such that $\nu\lrcorner\omega=0.$
Mind that \eqref{eq:absolutebc} might have a kernel, which is then given by the absolute de Rham cohomology $H_A^p(M)$ defined by 
$$H_A^p(M)=\{\phi\in \Omega^p(M)|\,\, d\phi=\delta_g\phi=0\,\, \text{on}\,\, M\,\, \text{and}\,\, \nu\lrcorner \phi=0\,\, \text{on}\,\, \partial M\}.$$ 
We also have a similar characterization for the first eigenvalue $\lambda_{1,p}^R$ of the Laplacian for the relative conditions with the corresponding relative cohomology $H_R^p(M)$ given by
$$H_R^p(M)=\{\phi\in \Omega^p(M)|\,\, d\phi=\delta_g\phi=0\,\, \text{on}\,\, M\,\, \text{and}\,\, \iota^*\phi=0 \,\, \text{on}\,\, \partial M\}.$$
By duality, the first eigenvalue for the boundary value problem \eqref{eq:absolutebc} on $p$-forms coincides with the first eigenvalue for the boundary value problem \eqref{eq:relativebc} on $(n-p)$-forms \cite{GS}. Another boundary problem of interest is the Dirichlet eigenvalue problem 
\begin{equation*}
\left\{\begin{array}{lll}\Delta \omega&=\lambda^D\omega &\textrm{ on }M\\ \omega &=0&\textrm{ on }\partial M.\end{array}\right.
\end{equation*}
For that problem, the first eigenvalue $\lambda^D_{1,p}$ -- which is necessarily positive by \cite[Thm. p. 445]{Anne89} -- is characterized by 
$$\lambda^D_{1,p}={\rm inf}\left\{\frac{\displaystyle\int_M |d\omega|^2d\mu_g}{\displaystyle\int_M|\omega|^2d\mu_g}\,,\,\, \omega\in\Omega^p(M)\setminus\{0\}\;\textrm{and}\;\omega_{|_{\partial M}}=0 \right\}.$$ 

One can easily see that when $\tau\to 0$, the Robin boundary problem (\ref{eq:BVPRobin}) reduces to the absolute boundary conditions. Also, when $\tau\to \infty,$ the problem (\ref{eq:BVPRobin}) reduces to the Dirichlet Laplacian. Now, we have the following bounds for the first eigenvalue $\lambda_{1,p}(\tau)$ of the Robin Laplacian on differential $p$-forms.

\begin{prop} \label{pro:neurodiri} Let $(M^n,g)$ be a compact Riemannian manifold with nonempty boundary.
Then, for any $\tau>0$ and all $p\in \{0,\cdots,n-1\}$, we have the double inequality
$$\lambda^N_{1,p}\leq\lambda_{1,p}(\tau)\leq\lambda^D_{1,p}.$$
\end{prop}  

{\noindent \it Proof.} The proof is based on the variational characterization of the first eigenvalue of each boundary value problem.
First, since $\tau>0$, for any $\omega\in\Omega^p(M)$ with $\nu\lrcorner\omega=0$, we have 
\[\int_M(|d\omega|^2 +|\delta_g\omega|^2)d\mu_g +\tau\int_{\partial M}|\iota^*\omega|^2 d\mu_g\geq\int_M(|d\omega|^2+|\delta_g\omega|^2)d\mu_g,\]
from which the left inequality follows.
Moreover, for any $\omega\in\Omega^p(M)$ with $\omega_{|_{\partial M}}=0$, we have 
\[\int_M(|d\omega|^2+|\delta_g\omega|^2)d\mu_g+\tau\int_{\partial M}|\iota^*\omega|^2d\mu_g=\int_M(|d\omega|^2+|\delta_g\omega|^2)d\mu_g\]
because of $\iota^*\omega=0$. Therefore $\lambda_{1,p}^D$ is the minimum of the same functional as that characterizing $\lambda_{1,p}(\tau)$ but taken on a smaller space (for $\omega_{|_{\partial M}}=0$ implies $\nu\lrcorner\omega=0$), which shows the right inequality.
\hfill$\square$

\noindent Next we establish a lower bound for the first eigenvalue of the Robin Laplacian on differential $p$-forms $\lambda_{1,p}(\tau)$ based on Theorem \ref{thm:bessel}. The lower term of the estimate depends on the so-called $p$-curvatures whose definition we recall. Let $\eta_1(x),\cdots,\eta_{n-1}(x)$ be the principal curvatures (i.e. eigenvalues of the Weingarten map $II$) at a point $x$ of the boundary $\partial M$ which can be assumed to satisfy $\eta_1(x)\leq\eta_2(x)\leq \cdots\leq\eta_{n-1}(x)$ up to reordering.
For any integer $p\in\{1,\cdots,n-1\},$ we define the $p$-curvatures $\sigma_p$ as $\sigma_p(x)=\eta_1(x)+\cdots+\eta_p(x).$ Clearly, one can check that for any two integer numbers $p$ and $q$ with $p\leq q$, we have that $\frac{\sigma_p(x)}{p}\leq \frac{\sigma_q(x)}{q}$ with equality if and only if $\eta_1(x)=\eta_2(x)=\cdots=\eta_q(x).$
From that remark follows the inequality $H\geq \frac{\sigma_p(x)}{p}$ for the mean curvature $H$ and for any $p\in\{1,\cdots,n-1\}.$ The Weingarten-endomorphism-field $II$ admits a canonical extension $II^{[p]}$ to $\Lambda^pT^*\partial M$ as follows: Given any $p$-form $\varphi$ on $\partial M,$ we define 
\begin{equation}\label{extensionweingarten}
(II^{[p]}\varphi)(X_1,\cdots,X_p)=\sum_{i=1}^p\varphi(X_1,\cdots,II(X_i),\cdots, X_p),
\end{equation}
where $X_i$ are vector fields on $\partial M$ for $i=1,\cdots,p.$ By a straightforward computation, it can be easily checked that the inequality
\begin{equation}\label{eq:weingarten}
\langle II^{[p]}\varphi,\varphi\rangle_x\geq \sigma_p(x)|\varphi|_x^2,
\end{equation}
holds pointwise. In the next theorem, we will denote by $\sigma_p$ as the infimum of $\sigma_p(x)$ over all $x\in\partial M.$ We have

\begin{thm}\label{estimatepform}
Let $(M^n,g)$ be a compact Riemannian manifold with smooth boundary.
Assume that $M$ has a nonnegative curvature operator and, for some $p\in \{1,\cdots,n-1\}$, the $p$-curvature of $\partial M$ is bounded from below by $\sigma_p>0$.
Fix any positive number $\tau_0<j_{\frac{n}{2}-1,1}$ and as before set $\alpha:=\tau_0\frac{J_{\frac{n}{2}}(\tau_0)}{J_{\frac{n}{2}-1}(\tau_0)}=\sum_{k\geq 1}\frac{2\tau_0^2}{j_{\frac{n}{2}-1,k}^2-\tau_0^2}.$
Then there exists an $\varepsilon>0$ such that, if $\tau>\sigma_p(\frac{\alpha}{2p}-1)-\varepsilon,$ we have 
$$\lambda_{1,p}(\tau)> \frac{\sigma_p^2}{2p^2}\tau_0^2.$$  
\end{thm}

{\noindent \it Proof.} We follow the idea of \cite[Thm. 3.3]{GS}. First of all by following the same steps as in \cite[Lemma 4.10]{GS}, it can be proved that, for any differential $p$-form $\omega$ satisfying the boundary conditions in \eqref{eq:BVPRobin}, we have 
$$\langle\nabla_\nu\omega,\omega\rangle=\langle II^{[p]}(\iota^*\omega),\iota^*\omega\rangle+\tau|\iota^*\omega|^2,$$
where $II^{[p]}$ is the canonical extension of the endomorphism $II$ defined previously.
Using the estimate \eqref{eq:weingarten}, we obtain
\begin{equation}\label{eq:integraldelta}
\int_M\Delta(|\omega|^2) d\mu_g=\int_{\partial M}\frac{\partial}{\partial\nu}(|\omega|^2) d\mu_g=2\int_{\partial M}\langle\nabla_\nu\omega,\omega\rangle d\mu_g\geq 2(\sigma_p+\tau)\int_{\partial M}|\omega|^2 d\mu_g.
\end{equation}
Recall now the Bochner formula for $p$-forms: $\Delta=\nabla^*\nabla+W_M^{[p]}$, where $W_M^{[p]}$ is the zero-order curvature term.
It is elementary to deduce from that formula the following scalar identity that is valid for any $p$-form $\omega$:
\begin{equation}\label{bochnerformula}
\langle\Delta\omega,\omega\rangle=|\nabla\omega|^2+\frac{1}{2}\Delta(|\omega|^2)+\langle W_M^{[p]}(\omega),\omega\rangle.
\end{equation}
Therefore, if $\omega$ is an eigenform associated to the eigenvalue $\lambda_{1,p}(\tau)$, then using the nonnegativity of $W_M^{[p]}$ (which is a consequence of that of the curvature operator of $M$), as well as $|\nabla\omega|^2,$ we obtain
\begin{equation}\label{ineeigenvalue}
\Delta(|\omega|^2)\leq 2\lambda_{1,p}(\tau) |\omega|^2.
\end{equation}  
Therefore, we are in the situation of Theorem \ref{thm:bessel} with $f=|\omega|^2.$ Recall here that the condition of the nonnegativity of the curvature operator on $M$ implies in particular that the Ricci curvature of $M$ is nonnegative. Also the mean curvature is bounded from below by $H_0:=\frac{\sigma_p}{p}>0.$ Let us from now on assume that $\sqrt{2\lambda_{1,p}(\tau)}< H_0 \tau_0.$ Then, we deduce from Theorem \ref{thm:bessel} that
$$\int_{\partial M} |\omega|^2 d\mu_g\geq  \sqrt{2\lambda_{1,p}(\tau)}\frac{J_{\frac{n}{2}-1}\left(\frac{\sqrt{2\lambda_{1,p}(\tau)}}{H_0}\right)}{J_{\frac{n}{2}}\left(\frac{\sqrt{2\lambda_{1,p}(\tau)}}{H_0}\right)}\int_M |\omega|^2 d\mu_g.$$

Note that $\displaystyle\int_{\partial M}|\omega|^2 d\mu_g$ cannot vanish since $\omega$ does not vanish identically by assumption.
But integrating \eqref{ineeigenvalue} over $M$ and using \eqref{eq:integraldelta}, we obtain  
$$\sqrt{2\lambda_{1,p}(\tau)}\geq 2(\sigma_p+\tau)\frac{J_{\frac{n}{2}-1}\left(\frac{\sqrt{2\lambda_{1,p}(\tau)}}{H_0}\right)}{J_{\frac{n}{2}}\left(\frac{\sqrt{2\lambda_{1,p}(\tau)}}{H_0}\right)}.$$
As we did in Corollary \ref{robinestimate}, we use the fact that the function $x\mapsto x\frac{J_{\nu+1}}{J_\nu}(x)$ is increasing on $]0,j_{\frac{n}{2}-1,1}[$ to obtain 
$$\sqrt{2\lambda_{1,p}(\tau)}\geq 2(\sigma_p+\tau)\frac{J_{\frac{n}{2}-1}\left(\frac{\sqrt{2\lambda_{1,p}(\tau)}}{H_0}\right)}{J_{\frac{n}{2}}\left(\frac{\sqrt{2\lambda_{1,p}(\tau)}}{H_0}\right)}> 2(\sigma_p+\tau)\frac{\sqrt{2\lambda_{1,p}(\tau)}}{\alpha H_0},$$
which implies $\sigma_p(\frac{\alpha}{2p}-1)> \tau.$
This contradicts the assumption on $\tau$.
Therefore $\frac{\sqrt{2\lambda_{1,p}(\tau)}}{H_0}\geq\tau_0$. The equality case of the latter estimate would provide the equalities in all the above inequalities. That means in particular the form $\omega$ has to be parallel, that is $\lambda_{1,p}(\tau)=0,$ which contradicts the fact that the Robin Laplacian has no kernel. Therefore $\frac{\sqrt{2\lambda_{1,p}(\tau)}}{H_0}>\tau_0$.
But since that inequality is strict, by continuity of $\tau\mapsto\lambda_{1,p}(\tau)$, there must exist an $\varepsilon>0$ such that, if $\tau>\sigma_p(\frac{\alpha}{2p}-1)-\varepsilon$, then $\frac{\sqrt{2\lambda_{1,p}(\tau)}}{H_0}>\tau_0$ keeps holding.
This concludes the proof of Theorem \ref{estimatepform}.
\hfill$\square$ 

\begin{remark}
{\rm By choosing $\tau_0=\nu-\frac{a_1}{2^\frac{1}{3}}\nu^\frac{1}{3}$ for $\nu=\frac{n}{2}-1$ with $n\geq 3$, one has that $\tau_0>\frac{n}{2}$ as in Remark \ref{remark:lowerbound}.
Therefore, for $\tau>\sigma_p(\frac{\alpha}{2p}-1),$ we have

\begin{equation}\label{estimationlambdasigma}
\lambda_{1,p}(\tau)>\frac{\sigma_p^2}{2p^2}\tau_0^2>\frac{n^2\sigma_p^2}{8p^2}.
\end{equation}

For Euclidean strictly $p$-convex domains (i.e. $\sigma_p>0$), we can easily see that the first eigenvalue $\lambda_{1,p}^N$ is positive.
This follows directly from integrating the Bochner formula (\ref{bochnerformula}) over $M$ and using  (\ref{eq:integraldelta}) for $\tau=0$, since a parallel form vanishing on the boundary must vanish identically. 
Combining the left inequality in Proposition \ref{pro:neurodiri} with the estimate $\lambda_{1,p}^N>\frac{\sigma_p^2}{8}$ established in \cite[p. 329]{GS} and which is valid for such domains, we get that $\lambda_{1,p}(\tau)>\frac{\sigma_p^2}{8}$ which is weaker than the estimate \eqref{estimationlambdasigma}.}
\end{remark}

\noindent In the following, we will use  the notation for the first eigenvalue of the Robin Laplacian on $0$-forms (which corresponds to functions) as 
$\lambda_{1,0}(\tau)=\lambda_1(\tau,M).$ We have 

\begin{cor} Let $(M^n,g)$ be a compact Riemannian manifold with smooth boundary. Assume that $M$ has a nonnegative curvature operator and the $p$-curvature of $\partial M$ is bounded by $\sigma_p>0$ for some $p\in \{1,\cdots,n-1\}.$ Then, we get the estimate 
$$\lambda_{1,p}(\tau)> \frac{\lambda_{1}(\tau,B_{H_0})}{2},$$
where $\lambda_{1}(\tau,B_{H_0})$ is the first eigenvalue of the scalar Robin Laplacian on the Euclidean ball whose boundary has mean curvature $H_0:=\frac{\sigma_p}{p}.$
\end{cor}

{\noindent \it Proof.} Fix any $\tau>0$. Recall that $0<\frac{\sqrt{\lambda_{1}(\tau,B_{H_0})}}{H_0}<j_{\frac{n}{2}-1,1}$, see Remark \ref{rem:eigenvalueball}.
Let $\tau_0:=\frac{\sqrt{\lambda_{1,0}(\tau,B_{H_0})}}{H_0}$.
Then $\alpha=\tau_0\frac{J_{\frac{n}{2}}(\tau_0)}{J_{\frac{n}{2}-1}}(\tau_0)=\tau_0\cdot\frac{\tau}{H_0\tau_0}=\frac{p\tau}{\sigma_p}$ by Corollary \ref{daners}. Therefore $\sigma_p(\frac{\alpha}{2p}-1)=\frac{\tau}{2}-\sigma_p$, in particular $\tau>\sigma_p(\frac{\alpha}{2p}-1)$.
It remains to apply Theorem \ref{estimatepform} to conclude the proof.
\hfill$\square$

 Next, we estimate the gap between the first eigenvalues of the Robin Laplacian for different orders when the manifold is isometrically immersed into the Euclidean space.
We mainly follow \cite[Thm 2.3]{GS}, \cite[Thm. 4]{RS1} and \cite{De}. 
In the following, for any given smooth unit normal vector field $\nu$ to $M$ in $\R^{n+m}$, we let $II_\nu$ be the associated Weingarten map, that is, it is the endomorphism field of $TM$ defined by 
\[\langle II_\nu(X),Y\rangle=\langle\nu,II(X,Y)\rangle\]
for all $X,Y$ tangent to $M$, where $II$ is the second fundamental form of the immersion.
We denote by $T^{[p]}$ the following endormorphism field of $\Lambda^pT^*\partial M$: given any local orthonormal basis $\{\nu_1,\ldots,\nu_m\}$ of $T^\perp M$, we let 
\[T^{[p]}:=\sum_{k=1}^m (II^{[p]}_{\nu_k})^2,\]
where each $II_{\nu_k}^{[p]}$ is the standard extension of the above $II_{\nu_k}$ to $\Lambda^pT^*\partial M$ as in \eqref{extensionweingarten}.
By \cite[Sec. 2]{GS}, the operator $T^{[p]}$ is well-defined, i.e independent on the chosen o.n.b. of $T^\perp M$, and furthermore self-adjoint as well as nonnegative.

\begin{thm}\label{t:estimgap}
Let $M^n\to \R^{n+m}$ be an isometric immersion where $(M^n,g)$ is a compact Riemannian manifold with $p$-convex boundary $\partial M$, that is $\sigma_p\geq 0$ for some $p\in \{1,\cdots,n-1\}.$
Then, for all $\tau>0,$ we have 
$$\lambda_{1,p}(\tau)-\lambda_{1,p-1}(\tau)\geq \frac{1}{p}\mathop{\rm inf}_M(W_M^{[p]}-T^{[p]}),$$ 
where $W_M^{[p]}$ is the Bochner operator in \eqref{bochnerformula}.
In particular, for a Euclidean $p$-convex domain $M\subset \R^{n+m}$, we have 
$$\lambda_{1,p}(\tau)\geq \lambda_{1,p-1}(\tau).$$
\end{thm}

{\noindent \it Proof.} Let $\omega$ be any eigenform of the Robin $p$-Laplacian associated with $\lambda_{1,p}(\tau)$.
Recall the boundary conditions $\nu\lrcorner\omega=0$ and $\nu\lrcorner d\omega=\tau\iota^*\omega$ valid for the form $\omega$ on $M$.
For each $i=1,\cdots, n+m,$ the unit parallel vector field $\partial_{x_i}$ on $\R^{n+m}$ splits into $\partial_{x_i}=(\partial_{x_i})^T+(\partial_{x_i})^\perp$ where $(\partial_{x_i})^T$ is the tangent part in $TM$ and $(\partial_{x_i})^\perp$ is the orthogonal one in $T^\perp M$.
Consider the $(p-1)$-form $(\partial_{x_i})^T\lrcorner \omega$ on $M$ which clearly satisfies $\nu\lrcorner ((\partial_{x_i})^T\lrcorner \omega)=0$ and apply to it the variational characterization in Proposition \ref{claim:varcharfirsteigenvRobin}. 
We get, for each $i$,  
\begin{equation}\label{eq:testpartialxi}\lambda_{1,p-1}(\tau)\int_M|(\partial_{x_i})^T\lrcorner \omega|^2 d\mu_g\leq \int_M(|d((\partial_{x_i})^T\lrcorner \omega)|^2 +|\delta_g((\partial_{x_i})^T\lrcorner \omega)|^2) d\mu_g+\tau\int_{\partial M}|(\partial_{x_i})^T\lrcorner \omega|^2 d\mu_g.\end{equation}
Now we want to sum over $i$.
We handle each term separately.
First, for any $p$-form $\alpha$ on $M$ and any local o.n.b. $\{e_1,\ldots,e_n\}$ of $TM$, we have 
\begin{eqnarray}\label{eq:sumpartialxi}
\nonumber\sum_{i=1}^{n+m}|\partial_{x_i}^T\lrcorner\alpha|^2&=&\sum_{i=1}^{n+m}\sum_{k,l=1}^n\langle\partial_{x_i}^T,e_k\rangle\langle\partial_{x_i}^T,e_l\rangle\langle e_k\lrcorner\alpha,e_l\lrcorner\alpha\rangle\\
\nonumber&=&\sum_{i=1}^{n+m}\sum_{k,l=1}^n\langle\partial_{x_i},e_k\rangle\langle\partial_{x_i},e_l\rangle\langle e_k\lrcorner\alpha,e_l\lrcorner\alpha\rangle\qquad\textrm{ since }\partial_{x_i}^\perp\perp TM\\
\nonumber&=&\sum_{k,l=1}^n\langle e_k,e_l\rangle \langle e_k\lrcorner\alpha,e_l\lrcorner\alpha\rangle\qquad\textrm{ since }(\partial_{x_i})_i\textrm{ is an o.n.b. of }\R^{n+m}\\
\nonumber&=&\sum_{k=1}^n\langle e_k\lrcorner\alpha, e_k\lrcorner\alpha\rangle\\
\nonumber&=&\langle\alpha,\sum_{k=1}^ne_k\wedge(e_k\lrcorner\alpha)\rangle\\
&=&p|\alpha|^2\qquad\textrm{ since }\sum_{k=1}^ne_k\wedge(e_k\lrcorner\alpha)=p\alpha.
\end{eqnarray}
Note that this remains valid pointwise, in particular also along $\partial M$.
To compute the term involving the exterior derivative, we make use of the Cartan formula and write, for each $i$,
\[d(\partial_{x_i}^T\lrcorner\omega)=L_{\partial_{x_i}^T}\omega-\partial_{x_i}^T\lrcorner d\omega.\]
By \cite[Eq. (4.3) p. 337]{GS}, we can split the Lie derivative as follows:
\begin{equation}\label{eq:splitLie}L_{\partial_{x_i}^T}\omega=\nabla_{\partial_{x_i}^T}\omega+II_{\partial_{x_i}^\perp}^{[p]}\omega.\end{equation}
As a first consequence, choosing $(e_k)_{1\leq k\leq n}$ and $(\nu_s)_{1\leq s\leq m}$ to be local o.n.b. of $TM$ and $T^\perp M$ respectively,
\begin{eqnarray*}
\sum_{i=1}^{n+m}|L_{\partial_{x_i}^T}\omega|^2&=&\sum_{i=1}^{n+m}|\nabla_{\partial_{x_i}^T}\omega|^2+|II_{\partial_{x_i}^\perp}^{[p]}\omega|^2+2\langle \nabla_{\partial_{x_i}^T}\omega,II_{\partial_{x_i}^\perp}^{[p]}\omega\rangle\\
&=&\sum_{i=1}^{n+m}\sum_{k,l=1}^n\langle\partial_{x_i}^T,e_k\rangle\langle\partial_{x_i}^T,e_l\rangle\langle\nabla_{e_k}\omega,\nabla_{e_l}\omega\rangle\\
&&+\sum_{i=1}^{n+m}\sum_{s,t=1}^m\langle\partial_{x_i}^\perp,\nu_s\rangle\langle\partial_{x_i}^\perp,\nu_t\rangle\langle II_{\nu_s}^{[p]}\omega,II_{\nu_t}^{[p]}\omega\rangle\\
&&+2\sum_{i=1}^{n+m}\sum_{k=1}^n\sum_{s=1}^m\langle \partial_{x_i}^T,e_k\rangle\langle\partial_{x_i}^\perp,\nu_s\rangle\langle\nabla_{e_k}\omega,II_{\nu_s}^{[p]}\omega\rangle
\end{eqnarray*}
As above, we can ignore the tangent and normal symbols in the sums on the r.h.s. and use the fact that $(\partial_{x_i})_i$ is an o.n.b. of $\R^{n+m}$:
\begin{eqnarray*}
\sum_{i=1}^{n+m}|L_{\partial_{x_i}^T}\omega|^2&=&\sum_{i=1}^{n+m}\sum_{k,l=1}^n\langle\partial_{x_i},e_k\rangle\langle\partial_{x_i},e_l\rangle\langle\nabla_{e_k}\omega,\nabla_{e_l}\omega\rangle\\
&&+\sum_{i=1}^{n+m}\sum_{s,t=1}^m\langle\partial_{x_i},\nu_s\rangle\langle\partial_{x_i},\nu_t\rangle\langle II_{\nu_s}^{[p]}\omega,II_{\nu_t}^{[p]}\omega\rangle\\
&&+2\sum_{i=1}^{n+m}\sum_{k=1}^n\sum_{s=1}^m\langle \partial_{x_i},e_k\rangle\langle\partial_{x_i},\nu_s\rangle\langle\nabla_{e_k}\omega,II_{\nu_s}^{[p]}\omega\rangle\\
&=&\sum_{k=1}^n|\nabla_{e_k}\omega|^2+\sum_{s=1}^m\langle II_{\nu_s}^{[p]}\omega,II_{\nu_s}^{[p]}\omega\rangle+2\sum_{k=1}^n\sum_{s=1}^m\underbrace{\langle e_k,\nu_s\rangle}_{0}\langle II_{\nu_s}^{[p]}\omega,II_{\nu_t}^{[p]}\omega\rangle\\
&=&|\nabla \omega|^2+\langle T^{[p]}\omega,\omega\rangle.
\end{eqnarray*}
Because of \eqref{eq:sumpartialxi}, we know that 
\[\sum_{i=1}^{n+m}|\partial_{x_i}^T\lrcorner d\omega|^2=(p+1)|d\omega|^2.\]
Moreover, still using the same local o.n.b.'s of $TM$ and $T^\perp M$ and decompositions of the $\partial_{x_i}$, we have 
\begin{eqnarray*}
\sum_{i=1}^{n+m}\langle L_{\partial_{x_i}^T}\omega,\partial_{x_i}^T\lrcorner d\omega\rangle&\bui{=}{\eqref{eq:splitLie}}&\sum_{i=1}^{n+m}\langle\nabla_{\partial_{x_i}^T}\omega,\partial_{x_i}^T\lrcorner d\omega\rangle+\langle II_{\partial_{x_i}^\perp}^{[p]}\omega,\partial_{x_i}^T\lrcorner d\omega\rangle\\
&=&\sum_{i=1}^{n+m}\langle\partial_{x_i}^T\wedge\nabla_{\partial_{x_i}^T}\omega,d\omega\rangle+\sum_{i=1}^{n+m}\sum_{k=1}^n\sum_{s=1}^m\langle\partial_{x_i}^\perp,\nu_s\rangle\langle\partial_{x_i}^T,e_k\rangle\langle II_{\nu_s}^{[p]}\omega,\nabla_{e_k}\omega\rangle\\
&=&|d\omega|^2+\sum_{k=1}^n\sum_{s=1}^m\underbrace{\langle\nu_s,e_k\rangle}_{0}\langle II_{\nu_s}^{[p]}\omega,\nabla_{e_k}\omega\rangle\\
&=&|d\omega|^2.
\end{eqnarray*}
It remains to notice that, since $\partial_{x_i}$ is parallel on $\R^{n+m}$, the covariant derivative $\nabla(\partial_{x_i}^T)$ is a symmetric endomorphism of $TM$ and therefore $\delta_g((\partial_{x_i})^T\lrcorner \omega)=-(\partial_{x_i})^T\lrcorner \delta_g\omega$.
Summing up \eqref{eq:sumpartialxi} over $i$, we obtain
\begin{eqnarray*}
\lambda_{1,p-1}(\tau)p\int_M|\omega|^2d\mu_g&\leq&\int_M\sum_{i=1}^{n+m}|L_{\partial_{x_i}^T}\omega|^2+|\partial_{x_i}^T\lrcorner d\omega|^2-2\langle L_{\partial_{x_i}^T}\omega,\partial_{x_i}^T\lrcorner d\omega\rangle d\mu_g\\
&&+(p-1)\int_M|\delta_g\omega|^2d\mu_g+\tau p\int_{\partial M}|\omega|^2d\mu_g\\
&\leq&\int_M|\nabla\omega|^2+\langle T^{[p]}\omega,\omega\rangle+(p+1)|d\omega|^2-2|d\omega|^2\\
&&+(p-1)\int_M|\delta_g\omega|^2d\mu_g+\tau p\int_{\partial M}|\omega|^2d\mu_g\\
&\leq &\int_M\left(|\nabla\omega|^2+\langle T^{[p]}\omega,\omega\rangle+(p-1)(|d\omega|^2+|\delta_g\omega|^2)\right)d\mu_g+\tau p\int_{\partial M}|\omega|^2d\mu_g.
\end{eqnarray*}
Applying formula \eqref{bochnerformula} and using \eqref{eq:integraldelta}, we deduce that 
\begin{eqnarray*}
\lambda_{1,p-1}(\tau)p\int_M|\omega|^2d\mu_g&\bui{\leq}{\eqref{bochnerformula}}&\int_M\left(\langle\Delta\omega,\omega\rangle-\frac{1}{2}\Delta(|\omega|^2)-\langle W_M^{[p]}\omega,\omega\rangle+\langle T^{[p]}\omega,\omega\rangle\right)d\mu_g\\
&&+(p-1)\int_M(|d\omega|^2+|\delta_g\omega|^2)d\mu_g+\tau p\int_{\partial M}|\omega|^2d\mu_g\\
&\bui{\leq}{\eqref{eq:integraldelta}}&\lambda_{1,p}(\tau)\int_M|\omega|^2d\mu_g-(\sigma_p+\tau)\int_{\partial M}|\omega|^2d\mu_g\\
&&+\int_M\langle (T^{[p]}-W_M^{[p]})\omega,\omega\rangle d\mu_g\\
&&+(p-1)\int_M(|d\omega|^2+|\delta_g\omega|^2))d\mu_g+\tau p\int_{\partial M}|\omega|^2d\mu_g\\
&\bui{\leq}{\eqref{eq:laplacianomega2}}&\lambda_{1,p}(\tau)\int_M|\omega|^2d\mu_g-(\sigma_p+\tau)\int_{\partial M}|\omega|^2d\mu_g\\
&&+\int_M\langle (T^{[p]}-W_M^{[p]})\omega,\omega\rangle d\mu_g\\
&&+(p-1)\left(\lambda_{1,p}(\tau)\int_M|\omega|^2d\mu_g-\tau\int_{\partial M}|\omega|^2d\mu_g\right)+\tau p\int_{\partial M}|\omega|^2d\mu_g\\
&\leq&p\lambda_{1,p}\int_M|\omega|^2d\mu_g-\sigma_p\int_{\partial M}|\omega|^2d\mu_g-\int_M\langle (W_M^{[p]}-T^{[p]})\omega,\omega\rangle d\mu_g.
\end{eqnarray*}
After regrouping the terms and getting rid of the boundary term using the condition $\sigma_p\geq 0$, we deduce the desired inequality. Finally, when $M$ is a Euclidean domain, both $W_M^{[p]}$ and $T^{[p]}$ vanish.
This ends the proof of Theorem \ref{t:estimgap}.
\hfill$\square$

\noindent Independently of the above results, we establish a Gallot-Meyer-type estimate for the first eigenvalue of the Robin Laplacian.
Recall that the Steklov (or Dirichlet-to-Neumann) operator on $p$-forms is the pseudo-differential operator $DN_p$ on $\Lambda^pT^*\partial M$ defined for any $p$-form $\omega$ on $\partial M$ by 
\[DN_p\omega:=-\nu\lrcorner d\hat{\omega},
\]
where $\hat{\omega}\in\Omega^p(M)$ is the unique $p$-form on $M$ such that $\Delta\hat{\omega}=0$ on $M$ with the boundary conditions $\iota^*\hat{\omega}=\omega$ as well as $\nu\lrcorner\hat{\omega}=0$ on $\partial M$, see \cite[Sec. 1.1]{RS1}.
It can be shown \cite[Thm. 11]{RS1} that the operator $DN_p$ is elliptic and essentially selfadjoint, in particular it has a discrete real spectrum consisting only of eigenvalues of finite multiplicities.
Its smallest eigenvalue $\nu_{1,p}$ can be characterized by  \cite[Thm. 11]{RS1}
$$\nu_{1,p}={\rm inf}\left\{\frac{\displaystyle\int_M \left(|d\omega|^2 +|\delta_g\omega|^2\right)d\mu_g}{\displaystyle\int_{\partial M}|\omega|^2d\mu_g}\right\},$$
where $\omega$ runs over all non-identically vanishing $p$-forms on $M$ such that $\nu\lrcorner\omega=0$.
We have 

\begin{thm}\label{t:GallotMeyerRobin}
Let $(M^n,g)$ be a compact Riemannian manifold with nonempty smooth boundary.
Assume that the curvature operator of $M$ is bounded from below by some $\gamma>0.$
Let $\tau\geq-\frac{c}{c-1}\cdot\sigma_p$ for some $p\in\{1,\cdots,n-1\}$ where $\sigma_p$ is the $p$-curvature of $\partial M$ and $c={\rm max}(p+1,n-p+1).$
Then
$$\lambda_{1,p}(\tau)\geq p(n-p)\frac{c}{c-1}\gamma.$$
If $\tau<-\frac{c}{c-1}\cdot\sigma_p$, we have 
$$\lambda_{1,p}(\tau)\geq \frac{p(n-p)(\nu_{1,p}+\tau)}{\frac{c-1}{c}\nu_{1,p}-\sigma_p}\gamma,$$
where $\nu_{1,p}$ is the first eigenvalue of the Steklov operator defined on differential $p$-forms.
\end{thm}
{\noindent \it Proof.} We follow mainly the proof as in the usual case.
{Let $\omega\in\Omega^p(M)$ be an eigenform associated to the eigenvalue $\lambda_{1,p}(\tau)$.
Combining the Bochner formula, the pointwise inequality $W_M^{[p]}\geq p(n-p)\gamma$ and estimate $|\nabla\omega|^2\geq\frac{1}{p+1}|d\omega|^2+\frac{1}{n-p+1}|\delta_g\omega|^2$ that is valid for any $p$-form $\omega$ \cite{GallotMeyer75}, we have:
\begin{eqnarray*}
\lambda_{1,p}(\tau)\int_M|\omega|^2 d\mu_g=\int_M\langle\Delta\omega,\omega\rangle d\mu_g&=&\int_M|\nabla\omega|^2 d\mu_g+\frac{1}{2}\int_M\Delta(|\omega|^2)d\mu_g+\int_M\langle W_M^{[p]}(\omega),\omega\rangle d\mu_g\\
&\geq&\int_M\frac{1}{p+1}|d\omega|^2d\mu_g+\frac{1}{n-p+1}|\delta_g\omega|^2 d\mu_g\\&&+\frac{1}{2}\int_M\Delta(|\omega|^2)d\mu_g
+p(n-p)\gamma\int_M|\omega|^2d\mu_g\\
&\bui{\geq}{\eqref{eq:integraldelta}}&\frac{1}{c}\int_M(|d\omega|^2+|\delta_g\omega|^2)d\mu_g+(\sigma_p+\tau)\int_{\partial M}|\omega|^2d\mu_g\\&&+p(n-p)\gamma\int_M|\omega|^2 d\mu_g\\ 
&\bui{=}{\rm\eqref{eq:laplacianomega2}} &\frac{1}{c}\int_M\langle\Delta\omega,\omega\rangle d\mu_g-\frac{\tau}{c}\int_{\partial M}|\omega|^2d\mu_g\\&&+(\sigma_p+\tau)\int_{\partial M}|\omega|^2d\mu_g+p(n-p)\gamma\int_M|\omega|^2d\mu_g.
\end{eqnarray*}
Thus, we deduce that 
\begin{equation}\label{eq:estimationbochner}
\lambda_{1,p}(\tau)\left(1-\frac{1}{c}\right)\int_M|\omega|^2d\mu_g\geq \left(\sigma_p+\tau-\frac{\tau}{c}\right)\int_{\partial M}|\omega|^2d\mu_g+p(n-p)\gamma\int_M|\omega|^2d\mu_g.
\end{equation}
Note that $\sigma_p+\tau-\frac{\tau}{c}\geq0$ if and only if $\tau\geq-\frac{c}{c-1}\sigma_p$.
This concludes the proof of the first part.
To prove the other part of the theorem, pick any eigenform $\omega$ of the Robin Laplacian associated to the eigenvalue $\lambda_{1,p}(\tau)$.
Then we get after using Equation \eqref{eq:laplacianomega2} that 
$$(\nu_{1,p}+\tau)\int_{\partial M}|\omega|^2d\mu_g\leq \lambda_{1,p}(\tau)\int_M|\omega|^2 d\mu_g.$$ 
Now if $\tau<-\frac{c}{c-1}\cdot\sigma_p$, we combine the last inequality with the one in \eqref{eq:estimationbochner} to get 
$$
\lambda_{1,p}(\tau)\left(1-\frac{1}{c}\right)\geq \left(\sigma_p+\tau-\frac{\tau}{c}\right)\frac{\lambda_{1,p}(\tau)}{\nu_{1,p}+\tau}+p(n-p)\gamma,$$ 
which is the desired estimate. Notice that the condition required on $\tau$ gives in particular that $\sigma_p<0$. }
\hfill$\square$

\section{Appendix} 

\noindent In this section, we review some basic facts on the Bessel functions of the first and second kind and their properties. For more details, we can refer to \cite{AS,B,W}.

\noindent The following differential equation, known as Bessel's equation 
$$x^2y''+xy'+(x^2-\nu^2)y=0$$ 
has the general solution 
$$y=AJ_\nu(x)+BY_\nu(x)$$
where $J_\nu,$ called Bessel function of first kind, is given by the series 
$$J_\nu(x)=\sum_{k=0}^\infty \frac{(-1)^k\left(\frac{x}{2}\right)^{\nu+2k}}{k!\Gamma(\nu+k+1)},$$
where $\Gamma$ is the gamma function. The Bessel function of the second kind $Y_\nu$ is related to the first kind by the formula 
$$Y_\nu(x)=\frac{J_\nu(x){\rm cos}(\nu\pi)-J_{-\nu}(x)}{{\rm sin}(\nu\pi)}.$$
It is given by the following expansion
\begin{eqnarray*}
Y_\nu(x)&=&\frac{2}{\pi}J_\nu(x)\left({\rm ln}\frac{x}{2}+\gamma\right)-\frac{1}{\pi}\sum_{k=0}^{\nu-1}\frac{(\nu-k-1)!}{k!}\left(\frac{x}{2}\right)^{2k-\nu}\\&&+\sum_{k=0}^\infty \frac{(-1)^{k-1}[(1+\frac{1}{2}+\cdots+\frac{1}{k})+(1+\frac{1}{2}+\cdots+\frac{1}{k+\nu})]}{k!(k+\nu)!}\left(\frac{x}{2}\right)^{2k+\nu},
\end{eqnarray*}
where $\gamma\simeq0.5772157$ is the Euler constant. 

\noindent For integer values of $\nu$, we take the limit $\nu\to n.$ In this particular case, $J_\nu$ and $J_{-\nu}$ are not linearly independent. Indeed, one has $J_{-\nu}(x)=(-1)^\nu J_\nu(x)$ and $Y_{-\nu}(x)=(-1)^\nu Y_\nu(x).$ 
For all real values of $\nu$, the Bessel functions can be expressed in terms of Bessel functions of lower orders by the formulas  
\begin{equation}\label{derivativebessel}
\left\{
\begin{array}{lll}
\medskip
J_{\nu+1}(x)=\frac{2\nu}{x}J_{\nu}(x)-J_{\nu-1}(x), \\\\
J'_{\nu}(x)=\frac{1}{2}(J_{\nu-1}(x)-J_{\nu+1}(x)),\\\\
J'_{\nu}(x)=J_{\nu-1}(x)-\frac{\nu}{x}J_{\nu}(x),\\\\
J'_{\nu}(x)=\frac{\nu}{x}J_{\nu}(x)-J_{\nu+1}(x).
\end{array}
\right.
\end{equation}  
The functions $Y_\nu$ satisfy the same equations above as $J_\nu.$ We also have the following identities, known as Lommel's formulas \cite[p. 46 and 77]{W}, which relate Bessel functions of different orders. Namely, 
\begin{equation}\label{relationbessel}
\left\{
\begin{array}{lll}
\medskip
J_{\nu-1}(x)J_{-\nu}(x)+J_{\nu}(x)J_{-\nu+1}(x)=2\frac{{\rm sin} \pi\nu}{\pi x}\\
Y_\nu(x)J_{\nu+1}(x)-Y_{\nu+1}(x)J_{\nu}(x)=\frac{2}{\pi x}.
\end{array}
\right.
\end{equation} 
A transformed version of the Bessel differential equation shows that the Bowman equation \cite[p.117]{B}
\begin{equation}\label{diffequbowm}
y''(x)-\frac{2\alpha-1}{x}y'+\left(\beta^2\gamma^2x^{2\gamma-2}+\frac{\alpha^2-m^2\gamma^2}{x^2}\right)y=0
\end{equation}
has the following solution 
$$y(r)=\left\{
\begin{array}{lll}
\medskip
x^\alpha\left(AJ_m(\beta x^\gamma)+BY_m(\beta x^\gamma)\right)&& \text{for integer}\,\, m \\
x^\alpha\left(AJ_m(\beta x^\gamma)+BJ_{-m}(\beta x^\gamma)\right)&& \text{for noninteger}\,\, m.
\end{array}
\right.
$$

\end{document}